\documentclass[A4paper,11pt,reqno]{article}
\usepackage{amsmath}
\usepackage{amssymb,latexsym}
\usepackage[english]{babel}
\usepackage[latin1]{inputenc} %lettres accentu\'ees
\usepackage[T1]{fontenc} %couper les mots qui contiennent des lettres accentu\'ees.
\usepackage{ae}
\usepackage{color,graphicx}
\usepackage{subfigure}
\vspace{\baselineskip}

\newdimen\AAdi%
\newbox\AAbo%
\def\AArm{\fam0 \rm}%
\def\AAk#1#2{\setbox\AAbo=\hbox{#2}\AAdi=\wd\AAbo\kern#1\AAdi{}}%
\def\AAr#1#2#3{\setbox\AAbo=\hbox{#2}\AAdi=\ht\AAbo\raise#1\AAdi\hbox{#3}}%

\newcommand{\un}{{\AArm 1\AAk{-.9}{l}l}}

\newtheorem{thm}{Theorem}[section]

\newtheorem{lem}[thm]{Lemma}
\newtheorem{prop}[thm]{Proposition}
\newtheorem{rem}[thm]{Remark}

\def\qedbox{$\rlap{$\sqcap$}\sqcup$}
\def\finpreuve{\nobreak\hfill\penalty250 \hbox{}
           \nobreak\hfill\qedbox\par\medskip}

\newcommand{\J}{J_t\,,\;t\,{\ge}\,0}
\newcommand{\Jk}{J^k_t\,,\;t\,{\ge}\,0}
\newcommand{\B}{B_t\,,\;t\,{\ge}\,0}
\newcommand{\X}{X_t\,,\;t\,{\ge}\,0}
\newcommand{\Xk}{X^k_t\,,\;t\,{\ge}\,0}

\renewcommand{\Re}{\mathrm{Re\:}}
\renewcommand{\Im}{\mathrm{Im\:}}
\newcommand{\noi}{\noindent}

\def\sectio#1{\setcounter{equation}{0}\section{#1}}
\def\eqp#1{\ensuremath\,#1}

\begin{document}

% Titre -----------------------------------------------------------

\title{\bf L\'{e}vy processes:\\ Hitting Time, Overshoot and Undershoot\\
II - Asymptotic Behaviour}
\author{{\small{\text{\bf Bernard ROYNETTE}$^{(1)}$}},\
{\small {\text{\bf  Pierre VALLOIS}$^{(1)}$}} {\small and} {\small
{\text{\bf Agn\`{e}s VOLPI }$^{(2)}$}}}
\date{Manuscript P661 submitted to SPA, October 2004}
% fin du titre ---------------------------------------------------

% adresse ------------------------------------------------------------
\maketitle 

{\small
\noindent (1)\,\, Universit\'e Henri Poincar\'e, Institut de
Math\'ematiques Elie Cartan, B.P. 239, F-54506 Vand\oe uvre-l\`es-Nancy Cedex\\
\noindent (2)\,\, ESSTIN, 2 rue Jean Lamour, Parc Robert Bentz,
54500 Vand\oe uvre-l\`es-Nancy, France.}

\noindent e-mail:\\
roynette@iecn.u-nancy.fr~~vallois@iecn.u-nancy.fr~~agnes.volpi@esstin.uhp-nancy.fr

\vspace{1ex}

\begin{abstract}
Let ($X_t, \; t\ge0$) be a L\'{e}vy process started at $0$, with
L\'{e}vy measure $\nu$  and $T_x$ the first hitting time of level
\mbox{$x>0$}~: \mbox{$\displaystyle T_x:=\inf{ \{ t\ge 0;\;  X_t>x
\}}$}. Let $F(\theta,\mu,\rho,.)$ be the joint Laplace transform
of $(T_x, K_x, L_x)$: $\displaystyle F(\theta,\mu,\rho,x):=\mathbb
E \left(e^{-\theta T_x-\mu K_x-\rho
L_x}\un_{\{T_x<+\infty\}}\right) \eqp{,} $ where \mbox{$\theta\geq
0$}, $\mu\geq 0,\;\rho\geq0,\;x\geq 0$,
  \mbox{$\displaystyle K_x:=X_{T_x}-x$}  and  \mbox{$\displaystyle L_x:=x-X_{T_{x^-}}$}.\\
If we assume that $\nu$ has finite exponential moments we exhibit an asymptotic expansion for $F(\theta,\mu,\rho,x)$, as $x\rightarrow +\infty$. A limit theorem involving a normalization of the triplet $(T_x,K_x,L_x)$ as $x\rightarrow +\infty$, may be deduced.\\
At last, if $\nu_{|_{\mathbb R_+}}$ has finite moment of fixed order, we  prove that the ruin probability $\mathbb P(T_x<+\infty)$ has  at most a polynomial decay.
\end{abstract}
\vspace{1ex}

\textbf{Keywords:}  L\'evy processes, ruin problem, hitting time, overshoot, undershoot, asymptotic estimates, functional equation, characteristic function, exponential and polynomial decay.\\

\textbf{AMS 2000 Subject classification:}
60E10, 60F05, 60G17, 60G40, 60G51,60J35,  60J65, 60J75, 60J80.

\vspace{1ex}

%=============================================================

\sectio{Introduction}

%==========================================================

\noindent {\bf 1.1}  Let ($\X$) be a L\'{e}vy process started at
$0$, with decomposition :
$$X_t=\sigma B_t-c_0t+J_t \quad ;\quad t\geq 0,$$
where $c_0\in\mathbb R$, $\sigma >0$, $(\B)$ is a standard
Brownian motion
started at $0$, $(\J)$ is a pure jump process, independent of $(\B)$.\\
Let us  denote by $\nu$ its L\'{e}vy measure, and  $T_x$ the first
hitting time of level $x>0$~:
\begin{equation}
T_x:=\inf\;\{t\geq0\;;\; X_t>x\}\eqp{.}
\end{equation}

\noi In this paper we essentially focus  on the asymptotic
behaviour of $T_x$ as $x\rightarrow \infty$. Since $\displaystyle
T_x=\inf\{t\geq 0 ; X_t/ \sigma> x/\sigma \}$, we can assume that
$\sigma =1$. Consequently, we have :
\begin{equation}\label{int1a}
    X_t= B_t-c_0t+J_t \quad ;\quad t\geq 0.
\end{equation}
 \noi If $(\X)$ is a Brownian motion with drift $-c_0$ (i.e.
$J_t=0$), then (cf. \cite{Karatszas91} p.$197$)~:
\begin{equation} \label{int1b}
P(T_x<+\infty)= \left\{
\begin{array}{lc}
  1 & \mbox{if} \ c_0 \leq 0 \\
  \displaystyle e^{-2c_0x} & \mbox{if} \ c_0>0. \\
\end{array}
\right.
\end{equation}
%\mathbb P(T_x<+\infty)=1\quad {\rm if}\quad c_0>0\quad {\rm and
%}\quad  \mathbb P(T_x<+\infty)=e^{-\gamma x}\quad {\rm if}\quad
%c_0<0\eqp{,}

%
%with $\gamma=-2c_0 >0$.

 \noi  In particular, if $c_0>0$, there exists
$\gamma>0$, $C,C'>0$ such that
\begin{eqnarray}
    &  \mathbb P(T_x<+\infty)\leq C'e^{-\gamma
    x}\label{proba-ruine-inf-expo}\\
 & \qquad \quad \mathbb P(T_x<+\infty)_{\stackrel{\sim}{x\rightarrow
 +\infty}}
Ce^{-\gamma x}\eqp{.} \label{proba-rui-equivalent}
\end{eqnarray}

%& \qquad \ \ \
%\begin{alignat}{2}
%\label{proba-ruine-inf-expo}
%&(i) &\quad & \mathbb P(T_x<+\infty)\leq Ce^{-\gamma x}\\
%&(ii) &\quad & \mathbb P(T_x<+\infty)\sim_{x\rightarrow +\infty}
%Ce^{-\gamma x}\eqp{.} \label{proba-rui-equivalent}
%\end{alignat}
%
\noi A lot of authors have generalized (\ref{proba-ruine-inf-expo})
when $(\X)$ is either a L\'{e}vy process, or a diffusion, see for instance
\cite{GerberLandry98},  \cite{FurrerSchmidli94},  \cite{Schmidli95},
\cite{Schmidli01},  \cite{DufresneGerber90}, \cite{KluppelbergKyprianou03}.\\
Suppose that $(\X)$ is a L\'{e}vy process with  no negative jump
(i.e. the support of its L\'{e}vy measure is included in
$\mathbb{R}_+$) and  the characteristic exponent $\psi$ of $(\X)$)
has a positive zero $\gamma_0$. Then,  in \cite{Doney91}, it has
proved that (\ref{proba-rui-equivalent}) holds, with
$\gamma=\gamma_0$ and  $\displaystyle
C=-\frac{\psi'(0)}{\psi'(\gamma_0)}$. This result was extended in
\cite{BertoinDoney94}, without any restriction
 on the support of $\nu$. However in this case  $C$ is not explicit.\\
 %[1ex]
%
% --------------------------------

\noi {\bf 1.2} Since $(X_t)$ may jump, if $T_x<\infty$, we have :
$X_{T_x}\geq x \geq X_{T_x-}$, where $X_{T_x-}$ denotes the
left-limit of $(X_t)$ at time $T_x$. Therefore it is interesting
to consider the overshoot $K_x$ and the undershoot $L_x$, these
r.v.'s  being defined on $\{T_x<+\infty\}$ as follows :
\begin{equation}
\label{def-de-undershoot-overshoot} K_x:=X_{T_x}-x\;;\qquad
L_x:=x-X_{T_x -}\eqp{.}
\end{equation}
\noi  One aim of this paper is a generalization of
(\ref{proba-rui-equivalent}). We first replace \mbox{$\mathbb
P(T_x\!<\!+\infty)$} by $F(\theta,\mu,\rho,x)$, where
$F(\theta,\mu,\rho,\cdot)$ is  the Laplace transform of the
triplet $(T_x,K_x,L_x)$~:
\begin{equation}
F(\theta,\mu,\rho,x):=\mathbb E \left(e^{-\theta T_x-\mu K_x -\rho
L_x} \un_{\{T_x<+\infty\}}\right)\eqp{.}
\end{equation}
%

%In particular $\displaystyle F(0,0,0,x)=\mathbb P(T_x<+\infty)$ is the ruin probability.\\
\noi Second,  we do not restrict ourselves  to  determine an
equivalent of $F(\theta,\mu,\rho,x)$, as $x$
  goes to infinity, we obtain (cf Theorem~\ref{thm:comp-asymp-de-F}) an asymptotic development of
this function. In our proof, the $x$-Laplace transform :
\begin{equation}
\widehat
F(\theta,\mu,\rho,q):=\int_0^{+\infty}e^{-qx}F(\theta,\mu,\rho,x)dx
\eqp{,}
\end{equation}
of $F(\theta,\mu,\rho,.)$ plays a central role.

% In a previous
%paper \cite{roynette-vallois-volpi-art1-2003}  we prove that
%$F(\theta,\mu,\rho,.)$ and  solves  some kind of integral
%equation, where $\widehat F(\theta,\mu,\rho,.)$ is the $x$-Laplace
%transform~:

\noi Let us explain the idea of our approach. For simplicity we
first
 consider \mbox{$\mathbb P(T_x<+\infty)$} (i.e. $\theta=\mu=\rho=0$).\\
 Obviously (\ref{proba-rui-equivalent}) is equivalent to
 $\displaystyle \lim_{x\rightarrow +\infty}e^{\gamma x}\mathbb P(T_x<+\infty)=C$. Suppose
  that condition holds, then
\begin{equation}
\label{sim-en-0} \int_0^{+\infty}e^{-q x}e^{\gamma x}\mathbb
P(T_x<+\infty)dx=\widehat
F(0,0,0,q-\gamma)_{\stackrel{\sim}{q\rightarrow0}}\frac{C}{q}\eqp{.}
\end{equation}
Hence $C$ is the residue of $\widehat F(0,0,0,q -\gamma)$ at $0$.

%\noi Conversely if (\ref{sim-en-0}) holds,
%(\ref{proba-rui-equivalent}) may be
% obtained  by a Tauberian theorem. Indeed in our context,
% the inversion of Mellin-Fourier is more convenient.

%$\gamma$ in (\ref{proba-rui-equivalent}),

\noi To deal with $F(\theta,\mu,\rho,x)$,  we have to suppose that
the jumps of $(X_t)$ are not too big :
\begin{equation}
\label{int-sur-R-} \int_{-\infty}^{-1}e^{-qy}\nu(dy)<+\infty\quad
\forall q>0\eqp{,}
\end{equation}
and
\begin{equation}\label{int11}
\int_{1}^{\infty}e^{sy}\nu(dy)<+\infty\quad \mbox{for some } \ s>0
\eqp{.}
\end{equation}
\noi With some additional assumptions (i.e.
(\ref{existence-de-kappa-bis}) and  (\ref{hyp-importante})), we
prove in Theorem~\ref{thm:comp-asymp-de-F}~
 :
\begin{align}
\label{dev-asymp} F(\theta,\mu,\rho,x) =&
 C_0(\theta,\mu,\rho) e^{- \gamma_0(\theta) x} +
       \sum_{i=1}^{p}a_i \left( C_i(\theta,\mu,\rho,x) e^{-\gamma_i(\theta) x}
+\overline{C_i}(\theta,\mu,\rho,x)
       e^{-\overline{\gamma_i}(\theta)x}\right)\nonumber\\
& +\mathrm{O}\left(e^{-Bx}\right)\eqp{,}
\end{align}
where $C_0(\theta,\mu,\rho)$ is a positive real number,
$C_1(\theta,\mu,\rho,x), \cdots, C_p(\theta,\mu,\rho,x)$ are $x$-polynomial
functions with values in $\mathbb C$,
$\left(\gamma_0(\theta),\gamma_1(\theta),
 \cdots,\gamma_p(\theta), \overline{\gamma_1}(\theta),
  \cdots, \overline{\gamma_p}(\theta)\right)$ are zeros of
   $\varphi-\theta$ (where $\varphi$ is the function defined by (\ref{def-phi}))
   and $\displaystyle a_i=\frac{1}{2}$ (resp.\ $1$)
   if  $\gamma_i(\theta)\in \mathbb R$  (otherwise).

\noi Few words about the proof of (\ref{dev-asymp}). Starting with
a functional equation verified by $\widehat F(\theta,\mu,\rho,.)$
(\cite{roynette-vallois-volpi-art1-2003},Theorem~$2.1$), we
determine the poles and the residues of $q \mapsto \widehat
F(\theta,\mu,\rho,q-\gamma_0(\theta))$. Then, using the
Mellin-Fourier inverse transformation we can recover $
F(\theta,\mu,\rho,.)$.

\noindent {\bf 1.3} In the same spirit, we investigate polynomial
decay of the probability of ruin. More precisely we prove in
Theorem~\ref{thm:majo-ruine} :
\begin{equation}\label{int12}
    P(T_x <\infty)\leq \frac{C}{1+x^n}, \quad \forall x\geq 0,
\end{equation}
where it is supposed that  $\displaystyle \int_0^\infty y^p
\nu(dy)<\infty$, for some $p\geq 2$, $n$ being the integer part of
$p-2$.

 \noi  Some results of this type have been obtained in  \cite{Paulsen01},
\cite{FrolovaKabanov01}, when
$(\X)$ is a diffusion .\\

\noi {\bf 1.4}  In section~\ref{sec:TCL} we give two applications
of Theorem~\ref{thm:comp-asymp-de-F}. The first one concerns
 the asymptotic behaviour of the triplet $(T_x,K_x,L_x)$ as $x\rightarrow\infty$.
Since $x\rightarrow T_x$ is non-decreasing, the first component
$T_x$ has to be normalized in $\widehat{T}_x$. We prove in
Theorem~\ref{thm:TCL-profit-net} (see also Remark~\ref{rem:TCL})
that $(\widehat{T}_x,K_x,L_x)$ converges in distribution, as
$x\rightarrow +\infty$. The first component $\widehat{T}_x$
converges to a Gaussian distribution. Morover  time and position
become asymptotically independent. In
\cite{KluppelbergKyprianou03}, the limit of distribution of the
overshoot has been investigated. A. Gut \cite{Gut88} has also
studied the convergence of the triplet, the L\'{e}vy process
$(\X)$ being replaced by a random walk.

\noi The asymptotic development of $F(\theta , \mu, \rho,x)$ given
in Theorem~\ref{thm:comp-asymp-de-F}  gives the rate of
convergence of $\widehat{T}_x$ to the Gaussian distribution (see
Remark \ref{rasymp3} for details).
%
%===========================================================

\sectio{Asymptotic expansion  of  $F(\theta,\mu,\rho,.)$}
\label{sec:comp-asymp-de-F}

%============================================================

We keep  the notation given  in the Introduction. Let $\psi$ the
characteristic exponent of $(\X)$, i.e. $\mathbb
E\left(e^{qX_t}\right)=e^{t\psi(q)}$. For our purpose it is more
convenient to deal with the function $\varphi$, where
$\varphi(q)=\psi(-q)$. L\'{e}vy Khintchine formula gives~:
\begin{equation}
\label{def-phi}
\varphi(q)=\frac{q^2}{2}+cq+\int_{\mathbb R}\left(e^{-qy}-1+ q y \un_{\{|y|<1 \}}\right)\nu(dy)\eqp{.}
\end{equation}
In this section, it is supposed that $\nu$ verifies :
\begin{equation}
\label{hyp-sur-r-bis}
r_\nu >0 \quad ({\rm may\; be\; }\quad r_\nu=+\infty)\eqp{,}
\end{equation}
\begin{equation} \label{cond-sur-r-nu-*}
r_\nu^* =\infty \eqp{,}
\end{equation}
\noi  where
\begin{equation}
\label{def-de-r-bis} r_\nu := \sup \big\{s \geq 0 ;\;
\int_{1}^{+\infty} e^{sy} \nu (dy )\;<\;+\infty \big\}\in
[0,+\infty]\eqp{,}
\end{equation}
\begin{equation}
\label{def-r-nu-*} r_\nu^*:=\sup\;\big\{q\in\mathbb
R\,;\;\int_{-\infty}^{-1}e^{-qy}\nu(dy)<+\infty\big\}\eqp{,}
\end{equation}
with the convention $\sup \emptyset =0$.

\noi Note that (\ref{cond-sur-r-nu-*}) is equivalent to
(\ref{int-sur-R-}).

\noi Consequently, in our setting,  $\varphi$ is defined on
$]-r_\nu,+\infty[$ and $\displaystyle
\int_{|y|>1}|y|\nu(dy)<\infty$. In particular :
\begin{equation}\label{asym1}
    E[|X_1|]<\infty.
\end{equation}
\noi Assumptions  (\ref{cond-sur-r-nu-*}) and
(\ref{hyp-sur-r-bis}) will be needed throughout this section,
therefore we do not repeat them in the statement of the results.

\noi  Let us briefly sketch  the proof of our main result
(Theorem~\ref{thm:comp-asymp-de-F}). We start with the functional
equation verified by $\widehat F(\theta,\mu,\rho,.)$. This result
was established in \cite{roynette-vallois-volpi-art1-2003},   and
the formula is recalled in (\ref{equa-Laplace}). To recover
$F(\theta,\mu,\rho,.)$ we use the Mellin-Fourier inverse
transformation. This leads us prove that some modification (namely
the function $\widehat{\widetilde{F}}(\theta,\mu,\rho,\cdot)$
defined by (\ref{asym3}))  of  $t\rightarrow \widehat
F(\theta,\mu,\rho,q+it)$  is an integrable function. We have also
  to finely analyze the zeros of $\varphi-\theta$.

 % \noi Before stating Theorem~\ref{thm:comp-asymp-de-F}, we need   some
 % recalls from \cite{roynette-vallois-volpi-art1-2003}
 %concerning the zeros of the function
 %$\varphi_\theta=\varphi-\theta$.
\noi As it is briefly explained in 1.2 of the Introduction, the
zeros of the function $\varphi_\theta=\varphi-\theta$ play a
important role.

 \noi Let us start with the negative zeros of $\varphi_\theta$.  We assume
 (cf.  Annex of  \cite{roynette-vallois-volpi-art1-2003},
   Figures $3.$a, $4.$a and $5$) that  there exists $\kappa >0$, such that
   for any $\theta \in[0,\kappa]$
   :
\begin{equation}
\label{existence-de-kappa-bis} \exists\; -\gamma_0(\theta)\in
]-r_\nu,0] \quad {\rm satisfying \; }\quad
\varphi(-\gamma_0(\theta))=\theta\eqp{.}
\end{equation}
More precisely~:
\begin{alignat}{2}
\label{existence-de-kappa-bis-point-i}
&(i)\quad  -\gamma_0(\theta)< 0,&& \quad  {\rm if}\quad \theta > 0\eqp{,} \\[1ex]
&\label{existence-de-kappa-bis-point-ii}
(ii)  \quad  -\gamma_0(0)<0,&&\quad {\rm if} \quad  \theta=0 \quad {\rm and}\quad  \mathbb E(X_1)<0\eqp{,}\\[1ex]
&(iii)
\label{existence-de-kappa-bis-point-iii}
 \quad -\gamma_0(0)=0, &&\quad {\rm if}  \quad \theta=0 \quad {\rm and}\quad  \mathbb E(X_1)\geq0
\eqp{.}
\end{alignat}
As for the positive zeros of $\varphi_\theta$, it is proved in  Annex of  \cite{roynette-vallois-volpi-art1-2003}~:
\begin{equation}
\label{existence-de-gamma0*-bis}
\forall\;\theta \geq 0,\qquad \exists\;\gamma_0^*(\theta)\geq 0\qquad {\rm such \; that}\qquad \varphi(\gamma_0^*(\theta))=\theta\eqp{.}
\end{equation}
More precisely~:
\begin{alignat}{2}
&(i)\label{existence-de-gamma0*-bis-point-i}
\quad \gamma_0^*(\theta)> 0 ,&& \qquad {\rm if}\quad \theta > 0 \eqp{,} \\[1ex]
&(ii)
\label{existence-de-gamma0*-bis-point-ii}
 \quad \gamma_0^*(0)>0, &&\qquad {\rm if}\quad \theta=0 \quad {\rm and}\quad \mathbb E(X_1)>0 \eqp{,}\\[1ex]
&(iii)
\label{existence-de-gamma0*-bis-point-iii}
 \quad \gamma_0^*(0)=0,&&\qquad {\rm if}  \quad \theta=0 \quad {\rm and}\quad   \mathbb E(X_1)\leq0
\eqp{.}
\end{alignat}
We set   $D_b:=\{q \in \mathbb C ;\quad {\rm such \; that }\quad
\Re q >-b\}$ and we introduce~:
\begin{equation}
B_\nu:=\sup{\{b > 0  \; ;\; \widehat {\nu_{|_{[1,+\infty[}}} \; {\rm admits \; a \; meromorphic\; extension\; to\; } \; D_b \}}\eqp{,}
\end{equation}
where \quad \mbox{$\displaystyle
\widehat {\nu_{|_{[1,+\infty[}}} (q): = \int_1^{+ \infty} e^{-qy}  \nu(dy)
$}.\\

\noindent Thanks to (\ref{hyp-sur-r-bis}), we note that $B_\nu\geq
r_\nu$ and $B_\nu$ may be equal to  $\infty$.  Then $\varphi$ has
a meromorphic extension to  $D_{B_\nu}$. For simplicity we still
denote  $\varphi$ this extension. Likewise,
 $\widehat {\nu_{|_{[1,+\infty[}}}$ is the meromorphic extension
   of $\widehat {\nu_{|_{[1,+\infty[}}}$ to  $D_{B_\nu}$.

\noindent Our last assumption on $\nu$ is~:
\begin{align}
\label{hyp-importante}
&\forall B\in]0,B_\nu[, \quad \exists K>0, \quad \exists R_0 >0, {\rm such \; that:}\nonumber\\
&|\widehat {\nu_{|_{[1,+\infty[}}} (q)|
%=\left| \int_1^{+ \infty} e^{-qy}  \nu(dy) \right|
 \leq K |q| \;,\quad {\rm for \;any}\; q \; {\rm such \;that}\quad \left\{
\begin{array}{c}
-B\leq \Re q\leq 0\\
|\Im q| \geq R_0
\end{array}\right.
\end{align}
A large class of measures $\nu$ satisfying previous hypothesis and
(\ref{hyp-importante}) will be given
 in  Remark~\ref{rem:comp-asymp} below.

\begin{rem} \begin{enumerate}
    \item Relations (\ref{hyp-sur-r-bis}) and (\ref{cond-sur-r-nu-*})
  allow us to replace $\displaystyle \widehat
{\nu_{|_{[1,+\infty[}}} (q)$ in (\ref{hyp-importante}), by  either
$\displaystyle \int_0^{+\infty}\!\!\!\!\!(e^{-qy}-1+ q y
\un_{\{|y|<1 \}})\nu(dy)$  or $\displaystyle
\int_{-\infty}^{+\infty}\!\!\!\!\!(e^{-qy}-1+ q y
 \un_{\{|y|<1 \}})\nu(dy)$.

    \item Since :
\begin{equation}
\label{nu-borne-sur-1+infini} \sup_{\Re q\geq 0}|\widehat
{\nu_{|_{[1,+\infty[}}} (q)|<+\infty\eqp{,}
\end{equation}
then (\ref{hyp-importante}) implies~:
\begin{equation}
\forall q \quad {\rm\; such\; that\; :}\quad \left\{
\begin{array}{c}
\Re q\geq -B\\
|\Im q| \geq R_0
\end{array}\right.\qquad
|\widehat {\nu_{|_{[1,+\infty[}}} (q)|
 \leq K |q| .
\end{equation}
\end{enumerate}

\end{rem}

\noi The key result concerning the complex  zeros of
$\varphi_\theta$
 is the following.
\begin{prop}
\label{prop:les-zeros-phi-theta}
We suppose that
%(\ref{hyp1}), (\ref{int-sur-R-}),  (\ref{hyp-sur-r-bis}),
(\ref{existence-de-kappa-bis}) and  (\ref{hyp-importante}) hold. Then  for any  $\theta \in [0,\kappa]$,  there exists  $\beta_\theta >0$  such that for any  $B \in ]0,\,B_\nu[$,  $\varphi_\theta$ admits a finite number of conjugated zeros in the strip \mbox{$\displaystyle
D_{B,\beta_\theta}:=\{q\in \mathbb C \;;\;  -B \leq  \Re q\leq \beta_\theta\}$}.\\
This set of zeros of $\varphi_\theta$ in $D_{B,\beta_\theta}$ is equal to~:
\begin{enumerate}
\item
$\displaystyle \{-\gamma_0(\theta), \; -\gamma_1(\theta), \; -\overline {\gamma_1}(\theta), \cdots, -\gamma_p(\theta), \;\displaystyle  -\overline {\gamma_p}(\theta)\}$ \quad if\quad
 $ \theta>0$\eqp{,}

\item
$\displaystyle \{0, \; -\gamma_0(0), \; -\gamma_1(0), \;\displaystyle -\overline {\gamma_1}(0), \cdots, -\gamma_p(0),\; \displaystyle -\overline {\gamma_p}(0)\}$
\quad if  \quad  $\theta=0$  and\; \mbox{$\mathbb E(X_1)<0$}\!\eqp{,}

\item
$\displaystyle \{-\gamma_0(0)=0, \; -\gamma_1(0), \; -\overline {\gamma_1}(0), \cdots, -\gamma_p(0), \;\displaystyle  -\overline {\gamma_p}(0)\}$ \; if
\;
 $\theta=0$ \;and \; \mbox{$\mathbb E(X_1)\geq 0$}\!\eqp{,}

where
\begin{equation}
\label{zeros-ordonnes}
-B < -\Re(\gamma_p(\theta))\leq \cdots\leq -\Re(\gamma_1(\theta))<
 -\gamma_0(\theta)\leq 0 \eqp{.} \\[2ex]
\end{equation}
\item  $-\gamma_0(\theta)$ is a simple (resp.\ double) zero of $\varphi_\theta$, if $\theta>0$ or $\theta=0$ and $\mathbb E(X_1)\neq0$ (resp.\ otherwise, i.e.\ $\theta=0$ and $\mathbb E(X_1)=0$).
\end{enumerate}
\end{prop}

\noindent {\bf Proof of Proposition~\ref{prop:les-zeros-phi-theta}}

\noindent $1.$\; It is clear that $\varphi_\theta$ is
holomorphic in \mbox{$\displaystyle \{q\in \mathbb C\;;\;\Re q>- r_\nu\}$}
and the only real zeros of $\varphi_\theta$ in this domain
are $-\gamma_0(\theta)$ and $\gamma_0^*(\theta)$.\\
\noindent $2.$\; We claim that $\varphi_\theta$ admits
only two  zeros, $-\gamma_0(\theta)$ and $\gamma_0^*(\theta)$,
in  \mbox{$\displaystyle \{q\in \mathbb C\;;\;-\gamma_0(\theta)
\leq Re q\leq \gamma_0^*(\theta)\}$}.\\
\begin{itemize}
 \item[a)] Suppose first that $-\gamma_0(\theta)<\Re q <\gamma_0^*(\theta)$. We have~:
\begin{equation*}
|e^{\varphi_\theta(q)}|=|\mathbb E(e^{-\theta X_1-\theta})|
\leq \mathbb E\left(e^{-\Re \theta X_1-\theta}\right)=e^{\varphi_\theta(\Re q)}<1\eqp{,}
\end{equation*}
since $\varphi_\theta<0$ on $]-\gamma_0(\theta),\gamma_0^*(\theta)[$.
 Then $\varphi_\theta(q)\neq 0$.\\
\item[b)] Let $q=-\gamma_0(\theta)+ib$, $b\in \mathbb R$. We compute $\varphi_\theta(q)$~:
\begin{align}
\varphi_\theta(q)=&\varphi_\theta(-\gamma_0(\theta))-\frac{b^2}{2}
+\int_{-\infty}^{+\infty} e^{\gamma_0(\theta)y}(cos (by)-1)\nu(dy)\nonumber\\
&+i\left(-b\gamma_0(\theta)+cb-\int_{-\infty}^{+\infty} e^{\gamma_0(\theta)y}
\left(sin (by)-by \un_{\{|y|<1 \}}\right)\nu(dy)\right)\eqp{.}
\end{align}
But $\varphi_\theta(-\gamma_0(\theta))=0$, then $\Re(\varphi_\theta(q))
\leq -\frac{b^2}{2}$. Consequently $\varphi_\theta(q)=0$ iff $b=0$.\\
\item[c)] The same reasoning applies to the case $q=-\gamma_0^*(\theta)+ib$.
\end{itemize}
\noindent $3.$\; Let us modify the decomposition of $\varphi_\theta$~:
\begin{equation}
\varphi_\theta(q)=\frac{q^2}{2}-cq-\theta+
\widehat{\nu_{|_{[1,+\infty[}}}(q)-\nu([1; +\infty[)
+\int_{-\infty}^1 \left(e^{-qy}-1+ q y \un_{\{|y|<1 \}}\right)\nu(dy)\eqp{.}
\end{equation}
Obviously (\ref{cond-sur-r-nu-*}) implies that
\mbox{$\displaystyle q\rightarrow \int_{-\infty}^1
\left(e^{-qy}-1+ q y \un_{\{|y|<1 \}}\right)\nu(dy)$} is
holomorphic in \goodbreak \noindent \mbox{$\displaystyle \{q\in \mathbb C\;;\;\Re q\leq
\beta\}$},
 for any $\beta\in \mathbb R$.
\begin{itemize}
\item[a)]
Suppose $\theta>0$ or $\theta=0$ and $\mathbb E(X_1)<0$. Then $\varphi_\theta$ admits $\gamma_0^*(\theta)$ as a unique positive zero.\\
The crucial point is the following~: assumption (\ref{hyp-importante}) implies that there exists $R>R_0>0$, $k>0$ such that \mbox{$\displaystyle |\varphi_\theta(q)|\geq kq^2$} for any $q$, $-B\leq \Re q\leq 0$, $\Im q>R$.\\
Proposition~\ref{prop:les-zeros-phi-theta} will be a direct
consequence of the two previous steps $1)$ and $2)$,
$\varphi_\theta(\overline{z})=\overline{\varphi_\theta(z)}$ and
$\varphi_\theta$ is meromorphic in $D_{B_\nu}$. \item[b)] The case
$\theta=0$ and $\mathbb E(X_1)\geq 0$, can be treated similarly.
\finpreuve
\end{itemize}

\begin{rem}
{\rm As it has been   pointed out in the previous proof,
assumption (\ref{hyp-importante}) may be replaced by~:
\begin{equation}
\exists \varepsilon>0\quad {\rm such\;that}\quad \inf_{q\in \mathcal D}\left|\frac{\frac{q^2}{2}-cq+\widehat{\nu_{|_{[1,+\infty[}}}(q)}{q^{1+\varepsilon}}\right|>0
\end{equation}
where  $\mathcal D=\{q\in \mathbb C\;;\; -B\leq \Re q\leq 0\;,\; |\Im q|\geq R_0\}$.
}
\end{rem}

% --------------------THEOREME--------------------------------------
\begin{thm}
\label{thm:comp-asymp-de-F} \quad We suppose
% (\ref{hyp1}), (\ref{int-sur-R-}), (\ref{hyp-sur-r-bis}) and
(\ref{existence-de-kappa-bis}) and (\ref{hyp-importante}). Let
 $B \in ]0,\,B_\nu[$, and  $\kappa>0$, small enough, such that
  (\ref{existence-de-kappa-bis}) holds  and  (\ref{zeros-ordonnes})
   is satisfied for any $\theta \in [0,\kappa]$ ( $p$ being  independent
   from  $\theta \in [0,\kappa]$).

\noindent Then for any $\theta \in[0,\kappa]$,
 $\mu\geq0$ and $\rho\geq 0$, there exists a positive number
  \mbox{$C_0(\theta,\mu,\rho)>0$} and complex $x$-polynomial functions
   $C_1(\theta,\mu,\rho,x), \cdots, C_p(\theta,\mu,\rho,x)$ such that

\noindent $F(\theta,\mu,\rho,x)$ has the following asymptotic expansion as $x \rightarrow +\infty$~:
\begin{align}\label{asymp4}
%\label{dev-asymp}
F(\theta,\mu,\rho,x) =&
 C_0(\theta,\mu,\rho) e^{- \gamma_0(\theta) x} +
       \sum_{i=1}^{p}a_i \left( C_i(\theta,\mu,\rho,x) e^{-\gamma_i(\theta) x}
+\overline{C_i}(\theta,\mu,\rho,x)
       e^{-\overline{\gamma_i}(\theta)x}\right)\nonumber\\
& +\mathrm{O}\left(e^{-Bx}\right)\eqp{,}
\end{align}
where $\displaystyle a_i=\frac{1}{2}$ if $\gamma_i(\theta)$ is
real and $a_i=1$ otherwise; the degree of $C_i(\theta,\mu,\rho,.)$
is $n_i-1$, where $n_i$ is the order of multiplicity of
$-\gamma_i(\theta)$ and $\mathrm{O}$  is uniform with respect to
$\mu\geq0$, $\rho\geq 0$  and $\theta \in [0,\kappa]$.
\end{thm}
%

% --------------------REMARQUE-----------------------
\begin{rem} \label{rem:comp-asymp}
\begin{enumerate}
    \item If we drop assumption (\ref{existence-de-kappa-bis}), then
 $\varphi_\theta$ has no zero located in $[-r_\nu, 0[$ and  the
 asymptotic expansion (\ref{asymp4}) reduces to~:
 %(\ref{dev-asymp})
\begin{equation}
F(\theta,\mu,\rho,x)=\mathrm{O} (e^{-(r_\nu-\varepsilon) x})\;,
\quad \varepsilon>0\eqp{.}
\end{equation}

    \item  Heuristically, no assumption on
the negative jumps is required to get :
\begin{equation}
F(\theta,\mu,\rho,x)\leq C e^{-\gamma_0(\theta) x}\eqp{.}
\end{equation}
However to obtain an equivalent, or an asymptotic development of
$F(\theta,\mu,\rho,x)$ when $x$ goes to infinity, it is natural to
suppose that the negative and the positive parts of the jumps of
$(X_t)$ are controlled. Our asymptotic development looks like a
perturbation theorem around the case of Brownian motion with
negative drift.

    \item We give three classes of measures $\nu$ satisfying (\ref{hyp-importante})~:
\begin{itemize}
 \item [a)]Suppose that  $\nu$ has finite exponential
moments~:
\begin{equation}
\label{moments-expo} \forall q \in \mathbb R \qquad
\int_{-\infty}^{+\infty}\left|e^{-qy}-1+ q y \un_{\{|y|<1
\}}\right|\nu(dy)<+\infty\eqp{.}
\end{equation}
In that case, $\widehat {\nu_{|_{[1,+\infty[}}}$  and  $\varphi$
are holomorphic functions in the whole plane $\mathbb C$, then
\mbox{$B_\nu=+\infty$}. Moreover, for any $B>0$~:
 \begin{equation}
\sup_{\Re q\geq-B}\left|\int_1^{+\infty} e^{-qy}\nu(dy)\right|\leq
\int_1^{+\infty} e^{B y}\nu(dy)<+\infty\eqp{.}
\end{equation}
Then (\ref{hyp-importante}) holds. Condition (\ref{moments-expo})
is realized if $\nu$ has a compact support.
%In particular if $\nu$ has compact support, the integral
%of  $|y|$  in a neighbourhood of  $0$ is finite then  $\nu$ verifies (\ref{moments-expo}).

\item [b)] Let  $\nu$ being a linear combination of  gamma
distributions~:
\begin{equation}
\nu(dy):=\sum_{i=1}^n \rho_i e^{-\beta_i
y}y^{m_i}\un_{\{y\geq0\}}dy\eqp{.}
\end{equation}
where   $\rho_i>0$, $\beta_i>0$ and $m_i\in \mathbb N$, for any
$i\in {\{1,2,\cdots,n\}}$,.

Since  the Laplace transform of $\nu$ is explicit, we obtain
immediately its  meromorphic extension to the whole plane
($B_\nu=+\infty$)  and (\ref{hyp-importante}). Indeed~:
\begin{equation}
\widehat \nu(q)=\int_0^{+\infty} e^{-qx} \sum_{i=1}^n \rho_i
e^{-\beta_i x}x^{m_i}dx=\sum_{i=1}^n \rho_i
\int_0^{+\infty}x^{m_i}e^{-(q+\beta i )x}dx \eqp{.}
\end{equation}
Setting  $y=(q+\beta_i)x$, we have~:
\begin{equation}
\widehat \nu(q)=\sum_{i=1}^n \frac{\rho_i}{(q+\beta_i)^{m_i+1}}
 \int_0^{+\infty}y^{m_i}e^{-y}dy
%=\sum_{i=1}^n \frac{\rho_i \Gamma(m_i+1)}{(q+\beta_i)^{m_i+1}}
=\sum_{i=1}^n \frac{\rho_i m_i!}{(q+\beta_i)^{m_i+1}}
  \eqp{.}
\end{equation}
This implies that $\widehat \nu$ and $\varphi$ are holomorphic in
$\mathbb C -{\{-\beta_1,\cdots,-\beta_n\}}$ and meromorphic on
$\mathbb C$.

\item [c)]Previous example may  be generalized taking~:
\begin{equation}
\nu(dy):=\phi(y) \un_{\{y\geq 0\}}dy\eqp{,}
\end{equation}
where $\phi \geq 0$, bounded on $[0,y_0]$, and for every $y\geq
y_0$~:
\begin{equation}
\phi(y):= \rho_0 e^{-\beta_0y}y^{m_0-1}+\sum_{i=1}^n (\rho_i
e^{-\beta_i y}+\overline {\rho_i} e^{-\overline{\beta_i}
y})y^{m_i-1} +\mathrm{O}(e^{-\beta_{n+1}\;y})\eqp{,}
\end{equation}
with  $y_0\geq 0$, $\rho_0 \geq 0$, $\beta_0 >0$,
$\Re(\beta_i)>0$, $\rho_i \in \mathbb C^*$, $m_i \in \mathbb N^*$
et \mbox{$\displaystyle \beta_{n+1} \geq\!\!\!
 \sup_{i \in {\{1,\cdots,n\}}}\!\!\!\!\! \Re (\beta_i)$}.
% and $\phi$ is a bounded function, constant on  $[y_0;+\infty[$ for some $y_0>0$.

\noindent Then  $\widehat \nu$ and   $\varphi$ are meromorphic
functions on \mbox{$\{q \;/\; \Re q\geq \beta_{n+1}\}$},
$B_\nu=\beta_{n+1}$ and (\ref{hyp-importante}) holds.

\item[d)]  If  $\nu=\nu_1 +\nu_2$,  $\nu_1$ and  $\nu_2$ verify
(\ref{hyp-importante}) then  $\nu$ also.
%\\[1ex]
\end{itemize}
\end{enumerate}
\end{rem}
% ------------------FIN DE LA REMARQUE------------------

\noindent {\bf Proof of Theorem~\ref{thm:comp-asymp-de-F}}

\noindent Recall (cf. Theorem~$2.1$ of
\cite{roynette-vallois-volpi-art1-2003}) that $\widehat
F(\theta,\mu,\rho,.)$ verifies~:
\begin{align}
\widehat {F}(\theta,\mu,\rho,q) = \frac{1}{\varphi (q)-\theta}&
\left( \frac{q-\gamma^*_0(\theta)} {2} +\int_0^{+
\infty}\left[\frac{e^{-(q+\rho)y}-e^{-\mu y}}{q+\rho-\mu} -
\frac{e^{-(\gamma^*_0(\theta)+\rho)y}-e^{-\mu
y}}{\gamma^*_0(\theta)+\rho-\mu} \right]
\nu(dy)\nonumber\right.\\[1ex]
& \left.+ RF(\theta,\mu,\rho,.)(q)- RF(\theta,\mu,\rho,.)
(\gamma^*_0(\theta)) \vphantom{\int}\right )\eqp{,}
\label{equa-Laplace}
\end{align}
where $\Re q>0$ and  $R$  denotes the operator~:
\begin{equation}
\label{def-de-R} Rh(q) : = \int_{-\infty}^0 \nu(dy)
\int_0^{-y}\left( e^{-q(b+y)}-1\right)  h(b) db  \eqp{.}
\end{equation}
\noindent For simplicity we only consider the case $\rho=0$. The
function  $F(\theta,\mu,0,\cdot)$ will be  written
$F(\theta,\mu,\cdot)$ for short.

\noi In the sequel,
 $B\in ]0,\,B_\nu[$ is supposed to be close to  $B_\nu$.

\paragraph{Step 1~: Replacing $F(\theta,\mu,.)$ by $\widetilde F(\theta,\mu,.)$
 and equation associated  with  $\widetilde F(\theta,\mu,.)$}\mbox{}

\noi We entend continuously  $F(\theta,\mu,.)$ to the whole line
in the following manner :
%$\widetilde F(\theta,\mu,.)$ dfinie et continue  sur $\mathbb R$ pa
\begin{equation}
\label{def-de-F-tild-chap7}
 \widetilde
F(\theta,\mu,x):=F(\theta,\mu,x)
\un_{[0;+\infty[}(x)+(1+x)\un_{[-1;+0]}(x) ,\qquad \forall x \in
\mathbb R \eqp{.}
\end{equation}
\noi Let  $\widehat{\widetilde F}(\theta,\mu.)$ be the Laplace
transform of  $\widetilde F(\theta,\mu,.)$~:
\begin{equation}\label{asym3}
\widehat{\widetilde F}(\theta,\mu,q):=
\int_{-\infty}^{+\infty}e^{-qx} \widetilde
F(\theta,\mu,x)dx\eqp{.}
\end{equation}
The advantage of using $\widetilde F(\theta,\mu,.)$  instead of
$F(\theta,\mu,.)$ lies in the fact that we shall prove in
Lemma~\ref{lem:Fourier-L1} below, that \mbox{$t\rightarrow
\widehat{\widetilde F}(\theta,\mu,q_1+it)$} is an integrable
function on $\mathbb R$, for any $q_1$ in $]0,\beta_\theta[$
($\beta_\theta$ being  defined in
Proposition~\ref{prop:les-zeros-phi-theta}).

\noindent Since \mbox{$\displaystyle
\int_{-1}^0(1+x)e^{-qx}dx=\frac{e^q-1-q}{q^2}$},  (\ref{equa-Laplace}) implies~:
\begin{align}
\label{FTC}
\widehat{\widetilde F}(\theta,\mu,q) = &
\frac{e^q-1-q}{q^2}\nonumber\\
&+
\frac{1}{\varphi (q)-\theta} \left (
\frac{q-\gamma^*_0(\theta)} {2} + \int_0^{+ \infty}\left[\frac{e^{-qy}-e^{-\mu y}}{q-\mu} - \frac{e^{-\gamma^*_0(\theta)y}-e^{-\mu y}}{\gamma^*_0(\theta)-\mu} \right] \nu(dy)
\nonumber\right.\\
&\left.\vphantom{\int}+RF(\theta,\mu,.)(q)-RF(\theta,\mu,.)(\gamma^*_0(\theta))\! \right)
\end{align}
We observe that
\begin{align}
&\frac{\varphi(q)-\theta}{q+\gamma^*_0(\theta)}=
\frac{\varphi(q)-\varphi(\gamma^*_0(\theta))}{q+\gamma^*_0(\theta)}\nonumber\\[1ex]
&=\frac{q-\gamma^*_0(\theta)}{2}+c \;\frac{q-\gamma^*_0(\theta)}{q+\gamma^*_0(\theta)}+\int_{-\infty}^{+\infty}\frac{e^{-qy}-
e^{-\gamma^*_0(\theta)y}+(q-\gamma^*_0(\theta))y \un_{\{|y|<1\}}}{q+\gamma^*_0(\theta)}\nu(dy)\eqp{,}
\end{align}
and
\begin{equation}
\frac{e^q-1-q}{q^2}+\frac{1}{q+\gamma^*_0(\theta)}=\frac{e^q-1}{q^2}-
\frac{\gamma^*_0(\theta)}{q(q+\gamma^*_0(\theta))}\eqp{.}
\end{equation}
Consequently for any $q$, such that $\Re q > 0$~:
\begin{align}
\label{F-tild-chapeau}
\widehat {\widetilde
F} (\theta,\mu,q)& = \frac{e^q-1}{q^2}-
\frac{\gamma^*_0(\theta)}{q(q+\gamma^*_0(\theta))}\nonumber\\[1ex]
&+
\frac{1}{\varphi(q)-\theta}
\left[-c \;\frac{q-\gamma^*_0(\theta)}{q+\gamma^*_0(\theta)}
-\int_{-\infty}^{+\infty}\frac{e^{-qy}-
e^{-\gamma^*_0(\theta)y}+(q-\gamma^*_0(\theta))y \un_{\{|y|<1\}}}{q+\gamma^*_0(\theta)}\nu(dy)
%\frac{1}{q+\gamma^*_0(\theta)}
%\int_{-\infty}^{+\infty}(e^{-qy}-e^{-\gamma^*_0(\theta)y})\nu(dy)
\right.\nonumber\\[1ex]
&+
\int_0^{+\infty}\left[\frac{e^{-qy}-e^{-\mu y}}{q-\mu}\nu(dy)
-\frac{e^{-\gamma^*_0(\theta)y}-e^{-\mu y}}{\gamma^*_0(\theta)-\mu}\right]\nu(dy)\nonumber\\[1ex]
&\left. +RF(\theta,\mu,.)(q)-RF(\theta,\mu,.)(\gamma^*_0(\theta))\vphantom{\int}\right]\eqp{.}
\end{align}

\noi It is clear that  (\ref{cond-sur-r-nu-*}) implies that  $RF(\theta,\mu,.)$  is an entire function on  $\mathbb C$ and all the integrals in (\ref{F-tild-chapeau}) are holomorphic at least in $D_{B_\nu}$. Consequently the right hand-side of (\ref{F-tild-chapeau}) is the meromorphic extension of  $\widehat{\widetilde F}(\theta,\mu,.)$ to  $D_{B_\nu}$. \\
Before ending step 1, we remark  that if  $\theta>0$ or
$\theta=0$ and  \mbox{$ \mathbb E(X_1)\neq 0$},
$\widehat{\widetilde F}(\theta,\mu,.)$ and
 $\widehat F(\theta,\mu,.)$ are  holomorphic in a neighbourhood of    $q=\gamma^*_0(\theta)$, although  $\gamma^*_0 (\theta)$ is a zero of $\varphi_\theta$, this value being a false singularity for  $\widehat{\widetilde F}(\theta,\mu,.)$.  But if  $\theta=0$ and  $\mathbb E(X_1)=0$,  then $\gamma^*_0(0)=0$ is a pole for  $\widehat {\widetilde F} (0,\mu,.)$.

\paragraph{Step 2~: $t\rightarrow \widehat{\widetilde F}(\theta,\mu,q_1+it)$ belongs to $\mathbb L^1(\mathbb R)$}\mbox{}\\

\noindent Let $\theta$ be a fixed element in $[0,\kappa]$,  $\beta_\theta>0$, given by  Proposition~\ref{prop:les-zeros-phi-theta} and   $q_1 \in]0,\beta_\theta[$.

\begin{lem}

\label{lem:Fourier-L1}

Under
%(\ref{hyp1}), (\ref{int-sur-R-}) and
(\ref{existence-de-kappa-bis}) and  (\ref{hyp-importante}),
%and let  $\;0<q_1<\beta_\theta$.
the function  $t\rightarrow \widehat{\widetilde F}(\theta,\mu,q_1+it)$ belongs to $\mathbb L^1(\mathbb R)$.
\end{lem}
To prove Lemma~\ref{lem:Fourier-L1}, we begin with stating few technical inequalities. These relations will be  also used  in the sequel.\\

\begin{lem}

\label{lem:integrales-bornees}
Let $\theta_1<\theta_2$.
We suppose
%(\ref{hyp1}), (\ref{int-sur-R-}) and
(\ref{hyp-importante}) is realized. Then~:
\begin{alignat}{2}
&(i) &&\quad \exists k>0\quad {\rm such \; that\; for\;any\;} q \; {\rm satisfying\;} \quad  \Re q \in [\theta_1,\theta_2], \quad {\rm we\;have~:}\nonumber\\[1ex]
&\label{prem-inegalite}{} &&\quad
\left| \int_{-1}^1 \left(e^{-qy}-1+q y \un_{\{|y|<1 \}}\right)\nu(dy)\right| \leq k|q|\eqp{,}\\[1ex]
&\label{autres-inegalites}{} &&\quad  \left| \int_{-1}^0 \left(e^{-qy}-1+q y \un_{\{|y|<1 \}}\right)\nu(dy)\right| \leq k|q|,\\[1ex]
& \label{des-inegalites}&&\quad
\left|\int_0^1 \left(e^{-qy}-1+q y \un_{\{|y|<1 \}}\right)\nu(dy)\right| \leq k|q| \eqp{,}\\[1ex]
&(ii) &&\quad
\forall d>0, \quad \exists k_0>0\quad {\rm such\; that} \quad \sup_{\Re q \leq d}
\left| \int_{-\infty }^{-1}e^{-qy}\nu(dy)\right| \leq k_0\\[1ex]
&(iii) &&\quad \forall A>0\quad \exists k_1>0 \quad {\rm
such\;that} \quad \forall q \quad {\rm satisfying\;:}\quad \left\{
\begin{array}{c}
-B\leq \Re q\leq A\\
|\Im q| \geq R_0
\end{array}\right.\eqp{,}\nonumber \\[1ex]
&
\label{majoration-de-nu} { } &&\quad
\left| \int_{-\infty}^{+\infty} \left(e^{-qy}-1+q y \un_{\{|y|<1 \}}\right)\nu(dy)\right|\; \leq\;k_1(1+|q|)\eqp{,}
\qquad \\[1ex]
&
\label{encadrement-de-phi}{ } &&\quad
\frac{1}{k_1}\;|q|^2 \;\leq\; |\varphi(q)|\; \leq\; k_1 \;|q|^2\qquad \\[1ex]
&(iv)
\label{majoration-de-RF} && \quad
\forall h \in \mathbb R, \quad
\sup_{\Re q\leq h}|RF(\theta,\mu,q)|<+\infty\eqp{.}
\end{alignat}
\end{lem}

\newpage

\noindent {\bf Proof of Lemma~\ref{lem:integrales-bornees}}

$(i)$  We set $q=a+ib$, where  $a \in [\theta_1,\theta_2]$. We have~:
\begin{align}
\left|\int_{-1}^1 \left(e^{-qy}-1+q y \un_{\{|y|<1 \}}\right)\nu(dy)\right|\leq&
\left|\int_{-1}^1 \left(e^{-(a+ib)y}-e^{-ay}+ i b y \un_{\{|y|<1 \}}\right)\nu(dy)\right|\nonumber \\
&+
\left|\int_{-1}^1 \left(e^{-ay}-1+ a y \un_{\{|y|<1 \}}\right)\nu(dy)\right|\eqp{.}
\end{align}

We deduce immediately (\ref{prem-inegalite}).

The two other inequalities (\ref{autres-inegalites}) and (\ref{des-inegalites})
 can be proved by the same way.\\[1ex]

$(ii)$ \quad  Let $q$ such that $\Re q \leq d$, then
\begin{equation}
\left|\int_{-\infty}^{-1} e^{-qy}\nu(dy)\right|
\leq \int_{-\infty}^{-1}e^{-\Re q y}\nu(dy)
\leq \int_{-\infty}^{-1}e^{-d y}\nu(dy)=k_0\eqp{.}
\end{equation}
Assumption (\ref{cond-sur-r-nu-*}) implies that $k_0$ is a finite constant.\\[1ex]

$(iii)$  \quad For any  $q$ satisfying  $-B \leq\Re q\leq A$ and
$|\Im q|\geq R_0$, we get

\begin{align}
\int_{-\infty}^{+\infty} \left(e^{-q y}-1+q y \un_{\{|y|<1 \}}\right)\nu(dy)&=
\int_{-\infty}^{-1}e^{-q y}\nu(dy)+\int_{-1}^{1}\left(e^{-q y}-1+q y
\un_{\{|y|<1 \}}\right)\nu(dy)\nonumber\\[1ex]
&+\int_{1}^{+\infty}e^{-q y}\nu(dy)
+\int_{-\infty}^{-1}\nu(dy)+\int_{1}^{+\infty}\!\!\!\nu(dy)\eqp{.}
\end{align}

Relation (\ref{majoration-de-nu}) follows from (\ref{hyp-importante}),  (\ref{prem-inegalite}) and (\ref{autres-inegalites}).

It is easy to check (\ref{encadrement-de-phi}) by virtue of the definition of $\varphi$ (cf. (\ref{def-phi})) and previous inequalities.

$(iv)$ Using the definition of $RF$ (cf. (\ref{def-de-R})), we deduce that $q\rightarrow RF(\theta,\mu,q)$ is an increasing function on $\mathbb R$ vanishing at $0$.

\noindent If  $q$ is a complex number such that $\Re q \leq h$, then
\begin{align}
|RF(\theta,\mu,q)| &\leq
\int_{-\infty}^0 \left( \int_0^{-y} e^{-\Re q (y+b)} F(\theta,\mu,b) db-y\right)\nu(dy) \nonumber\\[1ex]
& \leq k\; RF(\theta,\mu,h) \nonumber
\end{align}
This implies (\ref{majoration-de-RF}).
\finpreuve

\noindent {\bf Proof of Lemma~\ref{lem:Fourier-L1}}

\noi Proposition~\ref{prop:les-zeros-phi-theta} tells us
$\varphi_\theta$  has no zero in the strip
 \mbox{$\displaystyle \{q\in\mathbb C\;/\; 0<\Re q <\beta_\theta\}$}.  Hence   if $q_1 \in]0,\beta_\theta[$,  $\varphi(q_1+it)-\theta$ neither cancels and  (\ref{F-tild-chapeau}) implies that
\mbox{$\displaystyle t\rightarrow \widehat{\widetilde F}(\theta,\mu,q_1+it)$}
is a continuous function. Let us focus on  (\ref{F-tild-chapeau}).  Lemma~\ref{lem:integrales-bornees} implies  that all the numerators  in (\ref{F-tild-chapeau}) are bounded on the  line $\{q_1+it \; / \; t\in \mathbb R\}$,   and the denominators are less than $C|q^2|$ when $|q|\rightarrow +\infty$ where  $C>0$. This proves that
\mbox{$\displaystyle t\rightarrow \widehat{\widetilde F}(\theta,\mu,q_1+it)$} belongs to   $\mathbb L^1(\mathbb R)$.
\finpreuve

\paragraph{Step 3~: Proof of the asymptotic development (\ref{asymp4}),
 through the Mellin Fourier inverse transform}\mbox{}\\

\begin{figure} \centering
\scalebox{1} {\input{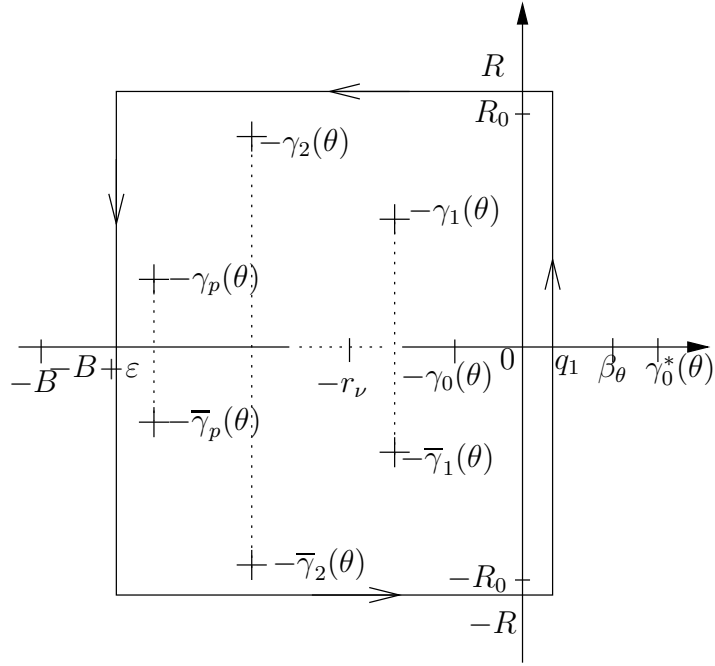}} \caption{The path
$\Gamma_{-B,q_1,R}$} \label{fig:chemin}
\end{figure}

\noindent  Proposition~\ref{prop:les-zeros-phi-theta} gives  the
existence of $\kappa>0$ and $B$ such that for any $\theta \in [0,\kappa]$,
 $\varphi_\theta$ does not vanish on \mbox{$\{-B+it \; / \;  t \in \mathbb R\}$}
 and $ B \neq \gamma^*_0(\theta)$.\\
Let $0<q_1<\beta_\theta$. Since
\mbox{$\displaystyle t \rightarrow \widehat{\widetilde F}
(\theta,\mu,q_1+it)$} belongs to  $\mathbb L^1(\mathbb R)$
(cf. Lemma~\ref{lem:Fourier-L1}), we are allowed to make use of the
Mellin Fourier inverse transform. So, for any $x\geq 0$~:
\begin{equation}
\label{eF-L1}
e^{-q_1 x}\widetilde F(\theta,\mu,x)=e^{-q_1 x} F(\theta,\mu,x)=
\frac{1}{2 \pi} \int_{-\infty}^{+\infty}e^{itx}
\widehat{\widetilde F}(\theta,\mu,q_1+it)
 dt\eqp{,}
\end{equation}
hence
\begin{equation}
 F(\theta,\mu,x)=
\frac{1}{2 \pi}
\int_{-\infty}^{+\infty}e^{(q_1+it)x} \widehat{\widetilde F}(\theta,\mu,q_1+it) dt
=-\;\frac{i}{2 \pi}\int_{\Gamma_{q_1}}e^{zx} \widehat{\widetilde F}(\theta,\mu,z)dz
\eqp{,}
\end{equation}
where  $\Gamma_{q_1}$ is the path~:
\begin{equation}
\Gamma_{q_1}:=\{z=q_1+it\quad {\rm such\; that}\quad t\in \mathbb R, \quad t\; {\rm increasing}\}\eqp{.}
\end{equation}
In   Proposition~\ref{prop:les-zeros-phi-theta}, it is proved there
exists \mbox{$R_1>R_0$}, such that  $\varphi_\theta$ has no zero in the two
 half-strips \mbox{$\displaystyle \{q\in\mathbb C\;/ -B\leq \Re q<\beta_\theta\;
 {\rm and }\; |\Im q|>R_1\}$}. In particular
 $\widehat{\widetilde F}(\theta,\mu,.)$ is holomorphic in this domain.\\
Let $\Gamma_{-B,q_1,R}$ be  the rectangular path (see
Figure~\ref{fig:chemin})~:

\begin{equation}
\Gamma_{-B,q_1,R}:=\Gamma_{q_1,R}\;\cup \; \Gamma_{R} \
;\cup \; \Gamma_{-B,R}\; \cup \; \Gamma_{-R}\eqp{,}
\end{equation}
where~:
\begin{alignat}{2}
&\Gamma_{q_1,R}&&:=\{q_1+it \; /\quad |t|\leq R,\quad t\;{\rm growing}\}\eqp{,} \\[1ex]
&\Gamma_{R}&&:=\{t+iR \; / \quad  -B\leq t\leq q_1,\quad t\;{\rm decreasing }\}\\[1ex]
&\Gamma_{-B,R}&&:=\{-B+it \; / \quad |t|\leq R,\quad t\;{\rm decreasing}\}\eqp{,} \\[1ex]
&\Gamma_{-R}&&:=\{t-iR \; / \quad  -B\leq t \leq q_1,\quad t\;{\rm growing}\}
\eqp{.}
\end{alignat}

%

%MODIFICATION FIGURE----------------------------------
%

\noindent Applying the residual theorem to the meromorphic extension
   of  $z \rightarrow e^{zx}\widehat{\widetilde F}(\theta,\mu,z)$ to
     $D_{B_\nu}$, we obtain, for any $R>R_1$~:
\begin{align}
\label{thm-residus}
&\int_{\Gamma_{-B,q_1,R}}e^{zx}\widehat{\widetilde F}(\theta,\mu,z)dz
=\nonumber\\
& 2i \pi \left[\vphantom{\sum_{i=1}^p} C_0(\theta,\mu)e^{-\gamma_0(\theta)x}\right.+
%\nonumber\\[1ex]
\left.\sum_{i=1}^p a_i \left(C_i(\theta,\mu,x)e^{-\gamma_i(\theta)x}+
\overline {C_i}(\theta,\mu,x) e^{-\overline {\gamma_i}(\theta))x}\right)\right]
\end{align}
where $\displaystyle a_i=\frac{1}{2}$ if  $-\gamma_i(\theta)$ is a real number  and  $a_i=1$ otherwise, and~:
\begin{align}
\label{def-constante-C0}
& C_0 (\theta,\mu) := Res
\left(\widehat{\widetilde F}(\theta,\mu,z);\; -\gamma_0(\theta)\right)\eqp{,}\\
& C_i (\theta,\mu,x) := e^{\gamma_i(\theta)x}Res
\left(e^{zx} \widehat{\widetilde F}(\theta,\mu,z);\; -\gamma_i(\theta)\right)
\label{def-des-constantes-Ci}
\eqp{,}
\end{align}
where  $Res(f(z);\gamma)$  denotes the residual of  $f$ at point  $\gamma$.
Let us remark (\ref{def-constante-C0}) is valid since $-\gamma_0(\theta)$
is a simple pole of $\widehat{\widetilde F}(\theta,\mu,.)$.\\
Since  \mbox{$\displaystyle z\rightarrow\frac{e^z-1-z}{z^2}$} has an holomorphic
 extension to the whole plan $\mathbb C$, identity (\ref{def-de-F-tild-chap7}) implies that~:
\begin{align}
\label{calcul-explicite-des-Ci}
& C_0(\theta,\mu)=Res
\left(\widehat F(\theta,\mu,z);\; -\gamma_0(\theta)\right)\eqp{,}\\
& C_i(\theta,\mu,x)= e^{\gamma_i(\theta)}
Res
\left(e^{zx}\widehat F(\theta,\mu,z);\; -\gamma_i(\theta)\right)\qquad
\forall i \in \{1,\cdots,p\}
\eqp{.}
\end{align}
We observe that  \mbox{$\displaystyle \overline {C_i}(\theta,\mu,x)=e^{\overline{\gamma_i}(\theta)x}Res\left(e^{zx}\widehat F(\theta,\mu,z);\; -\overline {\gamma_i}(\theta)\right)$}. We will prove in Remark~\ref{rem-ordre-multiplicity}  below that $x\rightarrow C_i(\theta,\mu,x)$ is a polynomial function.\\
\noindent Since \mbox{$z\rightarrow e^{zx}\widehat{\widetilde F}(\theta,\mu,z)$}
 belongs to $\mathbb L^1(\mathbb R)$ (cf. (\ref{eF-L1})), we have~:
\begin{align}
\label{F-contour}
 F(\theta,\mu,x)
&=-\;\frac{i}{2\pi}\; \lim_{R \rightarrow +\infty}\;
\int_{\Gamma_{q_1,R}}e^{zx}\widehat{\widetilde F}(\theta,\mu,z)dz
\nonumber\\[1ex]
&=-\;\frac{i}{2\pi}\; \lim_{R \rightarrow +\infty}\;\left[
\int_{\Gamma_{-B,q_1,R}}e^{zx}\widehat{\widetilde F}(\theta,\mu,z)dz
-\int_{\Gamma_{R}}e^{zx}\widehat{\widetilde F}(\theta,\mu,z)dz
\right.\nonumber\\[1ex]
&\left.\hphantom{=\frac{i}{2\pi}\:\;\quad}
-
\int_{\Gamma_{-B,R}}e^{zx}\widehat{\widetilde F}(\theta,\mu,z)dz
-\int_{\Gamma_{-R}}e^{zx}\widehat{\widetilde F}(\theta,\mu,z)dz \right]\eqp{.}
\end{align}
We claim that in the right hand-side of (\ref{F-contour}) the limits of the second term and the fourth one are null. As for the third limit, we have~:
\begin{equation}
\lim_{R \rightarrow +\infty}
\int_{\Gamma_{-B,R}} e^{zx} \widehat{\widetilde F}(\theta,\mu,z)dz=\mathrm{O}(e^{-Bx})\eqp{,}
\end{equation}
where $\mathrm{O}$ is uniform with respect to $\theta \in [0,\kappa]$ and
$\mu \in \mathbb R_+$.\\
Hence, as  $x\rightarrow+\infty$~:
\begin{align}
\label{partie-du-dev}
&F(\theta,\mu,x)=
C_0(\theta,\mu)e^{-\gamma_0(\theta)}
%\nonumber\\[1ex]&
+\sum_{i=1}^p\! a_i\! \left[\!
C_i(\theta,\mu,x)e^{-\gamma_i(\theta)x}
+\overline {C_i}(\theta,\mu,x) e^{ -\overline {\gamma_i}(\theta)x}
\right]
+\mathrm{O}(e^{-Bx})\eqp{.}
\end{align}

\noindent Let us determine the sign of $C_0(\theta,\mu)$. The asymptotic
 expansion (\ref{asymp4}) implies that~:
\begin{equation}
\lim_{x\rightarrow+\infty}e^{\gamma_0(\theta)x}F(\theta,\mu,x)=C_0(\theta,\mu)\geq 0
\eqp{.}
\end{equation}
The residual $C_0(\theta,\mu)$ at  $-\gamma_0(\theta)$  of $F(\theta,\mu,.)$ cannot vanish because $-\gamma_0(\theta)$ is a single pole of $\widehat F(\theta,\mu,.)$ when  $\theta>0$ or  $\theta=0$ and  $\mathbb E(X_1)\neq 0$.\\
If  $\theta=0$ and  $c=\mathbb E(J_1)$,  it is easy to see that  $C_0(0,\mu)\neq 0$ (cf. (\ref{res-0-mu-cas-3})).
\finpreuve

\begin{rem} \label{rem-ordre-multiplicity}
\begin{enumerate}
    \item We would like to show that
$C_i(\theta,\mu,\rho,.)$ is a polynomial function and determine
its maximal degree.  Assume that $-\gamma_i(\theta)$ is a zero of
$\varphi_\theta$ with  multiplicity $n_i$, and $\widehat
F(\theta,\mu,\rho,z)$ has the following asymptotic expansion in a
neighborhood of $-\gamma_i(\theta)$~:
\begin{equation}
\widehat
F(\theta,\mu,\rho,z)=\frac{K_{i,n_i}(\theta,\mu,\rho)}{(z+\gamma_i(\theta))^{n_i}}+
\frac{K_{i,n_i-1}(\theta,\mu,\rho)}{(z+\gamma_i(\theta))^{n_i-1}}+\cdots+\frac{K_{i,1}(\theta,\mu,\rho)}{z+\gamma_i(\theta)}+\cdots\eqp{.}
\end{equation}
Since
\begin{equation}
e^{zx}=e^{-\gamma_i(\theta)x}\left(1+(z+\gamma_i(\theta))x+\frac{(z+\gamma_i(\theta))^2}{2!}x^2+\cdots\right)
\end{equation}
then
\begin{equation}
C_i(\theta,\mu,\rho,x)=\frac{K_{i,n_i}(\theta,\mu,\rho)}{(n_i-1)!}x^{n_i-1}+\frac{K_{i,n_i-1}(\theta,\mu,\rho)}{(n_i-2)!}x^{n_i-2}+\cdots+K_{i,1}(\theta,\mu,\rho)\eqp{.}
\end{equation}

    \item \begin{itemize}
\item[(i)]  Suppose that   $-\gamma_i(\theta)$ is a single zero of
$\varphi_\theta$ (it is also a single pole of  $\widehat
F(\theta,\mu,.)$), then  $C_i(\theta,\mu,\rho,x)$ does  not depend
on $x$ and is given by the following~:
\begin{equation}
C_i(\theta,\mu,\rho)=Res \left( \widehat{ F}(\theta,\mu,\rho,z);\;
-\gamma_i(\theta)\right)\eqp{.}
\end{equation}
This situation  occurs if  $i=0$,  when  $\theta>0$ or  $\theta=0$ and  $\mathbb E(X_1)\neq 0$.\\
If the real part of  $-\gamma_i(\theta)$ is bigger than  $-r_\nu$, $C_i(\theta,\mu,\rho)$ can be determined as follows~: %through (\ref{equa-Laplace})~:
\begin{align}
\label{residu} C_i(\theta,\mu,\rho)=&
\frac{1}{\varphi'(-\gamma_i(\theta))}
\left[\frac{-\gamma_i(\theta)-\gamma^*_0(\theta)}{2}\right.\nonumber\\
&+\int_0^{+\infty}\left[\frac{e^{(\gamma_i(\theta)-\rho)y}-e^{-\mu
y}}{-\gamma_i(\theta)+\rho-\mu}
-\frac{e^{-(\gamma^*_0(\theta)+\rho)y}-e^{-\mu y}}{\gamma^*_0(\theta)+\rho-\mu}\right]\nu(dy)\nonumber\\[1ex]
&\left.
+RF(\theta,\mu,\rho,.)(-\gamma_i(\theta))-RF(\theta,\mu,\rho,.)(\gamma^*_0(\theta))
\vphantom{\int}\right] \eqp{,}
\end{align}
where it is supposed that $\displaystyle \frac{e^{ay}-1}{a}=y$ if $a=0$.\\
When $\Re(-\gamma_i(\theta))<-r_\nu$, the previous $\nu$-integrals and $\varphi$, have to be replaced by their meromorphic extensions.\\
 Recall that if $\mathbb E(X_1)<0$ (resp. $\mathbb E(X_1)>0$) then $\gamma^*_0(0)=0$  (resp. \mbox{$\gamma_0(0)=0$}).\\
\item[(ii)] If $\mathbb E(X_1)=0$, then $\gamma_0(0)=0$ is a
double zero of $\varphi$, but a simple pole of $\widehat
F(\theta,\mu,\rho,.)$. Thus, by a direct calculation, we have~:
\begin{align}
\label{res-0-mu-cas-3} C_0(0,\mu,\rho)&
=\frac{1}{\varphi''(0)}\left(1
-\frac{2}{(\rho-\mu)^2}\int_0^{+\infty}e^{-\rho
y}\left(1-e^{(\rho-\mu) y}+(\rho-\mu) y\right) \nu(dy)\right.
\nonumber\\[1ex]
& \hphantom{==}\left. -  2 \int_{-\infty}^0 \nu(dy)
\int_0^{-y}(y+b)F(0,\mu,\rho,b)db \right) \eqp{.}
\end{align}
In particular~:
\begin{align}
\label{res-0-cas-3} C_0(0,0,0)& = \frac{1}{\varphi''(0)}\left(1
+\int_0^{+\infty}y^2\; \nu(dy)
 - 2 \int_{-\infty}^0 \nu(dy) \int_0^{-y}\!\!\!\!(y+b)F(0,0,b)db
\right)\nonumber\\
&=1-\frac{2}{\varphi''(0)}\int_{-\infty}^0 \nu(dy)
\int_0^{-y}(y+b)F(0,0,0,b)db \eqp{.}
\end{align}
\end{itemize}
    \item
In \cite{BertoinDoney94}, it was proved that $\displaystyle
\mathbb P(T_x<+\infty)\sim C_0 e^{-\gamma_0 x}$ as $x\rightarrow
+\infty$, $C_0$ being indeterminate. We go a little bit further
since $C_0=C_0(0,0,0)$.
 \noindent If we suppose moreover that the support of $\nu$ is included in $[0,+\infty[$,
 then $RF(\theta,\mu,\rho,.)=0$ and

\begin{align} \label{calcul-des-coef}
& C_i(0,0,0)=-\; \frac{\varphi'(0)}{\varphi'(-\gamma_i(0))}, \quad {\rm if\;} \mathbb E(X_1)<0 \eqp{,}\\
&  C_i(0,0,0)=0 ,\; i\neq 0\;\;{\rm and }\;\; C_0(0,0,0)=1, \quad
{\rm if\;} \mathbb E(X_1)\geq 0 \eqp{.}
\end{align}

We recover the result given in \cite{Doney91}, i.e.\
(\ref{calcul-des-coef}) with $i=0$.
\end{enumerate}

\end{rem}

%==========================================================================

\sectio{A polynomial upper bound for the ruin probability}
\label{sec:majo-ruine}

%=========================================================================

In this section we only consider the ruin probability. For simplicity we note~:
\begin{equation}
F(x):=F(0,0,0,x)=\mathbb P(T_x<+\infty)\eqp{,}
\end{equation}
and its Laplace transform~:
\begin{equation}
\widehat F(q):=\int_0^{+\infty} e^{-q x}F(x) dx \qquad \forall q\in \mathbb C
\; ,  \; \Re(q)>0
\eqp{.}
\end{equation}
\noi We suppose in this section~:
\begin{equation}
\label{2-hypos}
\int_{|y|>1}|y|\nu(dy)<+\infty\qquad {\rm and}\qquad \mathbb E(X_1)=\int_{|y|>1}y\nu(dy)-c<0 \eqp{.}
\end{equation}
\noi Therefore  $\displaystyle \lim_{x\rightarrow +\infty}\mathbb
P(T_x<+\infty)=0$.

\noi  In   Section~\ref{sec:comp-asymp-de-F}  we have proved that,
under suitable assumptions, the ruin probability goes to $0$, with
exponential rate. Here the aim is to prove that under suitable
assumptions, $F(x)$ has a polynomial type rate of decay, as
$x\rightarrow \infty$.

\begin{thm}
\label{thm:majo-ruine}

 \noindent We suppose (\ref{2-hypos}) and
\begin{equation}
\label{moment-p-fini} \int_0^{+\infty} y^p \nu(dy)
<+\infty\;,\quad {\rm for\;some\;}\quad p\geq 2 \eqp{.}
\end{equation}
Let $n$ be the integer part of  $ p-2$, then
\begin{equation}
\label{majo-poly-de-F}
\forall x \in \mathbb R_+ \qquad  F(x) \leq \frac{C_n}{1+x^n}
\eqp{,}
\end{equation}
where $\displaystyle C_n>0$.
\end{thm}

\noindent {\bf Proof of Theorem~\ref{thm:majo-ruine}}\mbox{}

\noindent The proof will be divided into five parts.

%Nous prouvons que cette transforme de Fourier est $p-2$ fois
%drivable et que toutes ses drives jusqu' l'ordre $p-2$ sont dans
%$\mathbb L^1(\mathbb R)$.
%Nous concluons en utilisant alors le Thorme de Riemann-Lebesgue.\\[1ex]

\paragraph{{\bf Step 1}}\mbox{}
 We prove that  it is sufficient to consider a L\'{e}vy measure
  $\nu$ with support included in $[-k,+\infty[$, for some (finite) $k\geq0$.\\

\noi Assumption (\ref{2-hypos}) implies there exists
$k>0$ such that
\begin{equation}
\int_{-k}^{+\infty} \un_{\{|y|>1  \}} y \nu(dy)<c\eqp{.}
\end{equation}
Let  ($\Xk$) be a L\'{e}vy process with decomposition~:
\begin{equation}
X^k_t:=B_t-ct+J^k_t \qquad \forall t \geq 0
\eqp{,}
\end{equation}
where  $(\Jk$) is a  pure jump process with L\'{e}vy measure
\mbox{$\displaystyle \nu_k=\nu_{|_{[-k;+\infty[}}$}, independent
of $(\B)$. Moreover  $(\Jk$),  ($\J$) can be defined on the same
probability space, and  $\displaystyle J_t\leq J_t^k$, $\forall
t\geq 0$. Then  a.e. $\displaystyle X_t \leq X^k_t,\;\; \forall t
\geq 0$ and $\displaystyle T^k_x \leq T_x,\;\;  \forall x \geq 0 $
where
\begin{equation}
T^k_x:=\inf{\{t\geq 0\; / \; X^k_t>x\}}
\eqp{.}
\end{equation}
As a result
\begin{equation}
 \forall x\geq 0\qquad F(x)\leq F^k(x):=\mathbb P(T^k_x<+\infty)
\eqp{.}
\end{equation}
\noi This proves the claim.

\vskip5pt \noi  In the sequel  we  suppose that  the support of
$\nu$ is included in $[-k,+\infty[$, for some $k>0$.

%\begin{equation}
%\label{moment-2}
%\int_0^{+\infty}y^2\nu(dy)<+\infty
%\eqp{.}\\[2ex]
%\end{equation}

\paragraph{{\bf Step 2}}   $F$ belongs to $\mathbb L^1(\mathbb R_+)$.

\vskip4pt \noi We first prove~:
\begin{equation}
\label{sup-F-chapeau}
\sup_{0< q\leq q_0} |\widehat F(q)|<+\infty\;,\quad {\rm for\;some\;} q_0>0\eqp{.}\\[1ex]
\end{equation}
Choosing $\theta=\mu=\rho=0$ in (\ref{equa-Laplace}), we obtain~:
\begin{equation}
\label{trans-laplace-ruine}
\widehat F(q)=
 \frac{1}{\varphi(q)}
\left(\frac {q}{2}+\frac{1}{q} \int_0^{+\infty}(e^{-qy}-1+qy)\nu (dy)
+
RF(q) \right)\eqp{.}
\end{equation}
Let us determine the asymptotic  behaviour of the numerator and
the denominator, $q\rightarrow 0$~, we have :
\begin{align}
\label{equi-0-phi}
&\varphi(q) \sim_{q \rightarrow 0} q\varphi'(0) =q \left(c-\int_{|y|>1}y\nu(dy)\right),\\[1ex]
&\frac {q}{2}+\frac{1}{q} \int_0^{+\infty}(e^{-qy}-1+qy)\nu (dy)\sim_{q\rightarrow 0}
\frac{q}{2}\left[1+\int_0^{+\infty}y^2 \nu(dy) \right]
\label{equi-0-int}
\eqp{.}
\end{align}
Since the support of $\nu$ is included in $[-k,+\infty[$, $RF(q)$
can be simplified~:
\begin{equation}
RF(q)=\int_{-k}^0 e^{-qy}\nu(dy)\left[\int_0^{-y}e^{-q b}F(b)db-
\int_0^{-y}F(b)db\right]
\eqp{.}
\end{equation}
But $RF(0)=0$ and the derivative of $RF(q)$ is bounded, then
$\displaystyle |RF(q)|\leq  C q$, for any  \mbox{$\displaystyle 0\leq q \leq q_0$}.
(\ref{sup-F-chapeau}) follows immediately.\\[1ex]
It is now easy  to check that $F$ is in $\mathbb L^1(\mathbb
R_+)$. $F$ being positive~:
\begin{equation}
\int_0^{+\infty} F(x)dx
=\lim_{q\rightarrow 0}\int_0^{+\infty}e^{-qx}F(x)dx
\leq \sup_{0< q \leq q_0}\widehat F(q)<+\infty
\eqp{.}
\end{equation}
\vskip6pt\noi The function $F$ can be extended to the whole line,
setting $F(x)=0$,  for any   $ x \leq 0$. However $F$ may have a
jump at $0$. Let $\widetilde F$ be the following continuous
extension of $F$ :
\begin{equation}
\widetilde F(x):=F(x) \un_{[0;+\infty[}(x)+(1+x)\un_{[-1;+0]}(x),
\qquad \forall x \in \mathbb R \eqp{.}
\end{equation}
Let  $q\rightarrow \widehat{\widetilde F}(iq)$ be the Fourier transform of $\widetilde F$~:
\begin{equation}
\widehat{\widetilde  F}(iq)
=\int_{-\infty}^{+\infty}e^{-iqx}\widetilde F(x)dx
=\int_{-1}^{+\infty}e^{-iqx}\widetilde F(x)dx   \qquad \forall q\in \mathbb R\eqp{.}
\end{equation}

\paragraph{{\bf Step 3}}  $q\rightarrow \widehat{\widetilde F}(iq)$ is in  $ \mathbb L^1(\mathbb R)$.

\noi For simplicity, we suppose moreover that  $\displaystyle \int_{-1}^1|y| \nu(dy)< + \infty$.\\
\noindent Since $F\in  \mathbb L^1(\mathbb R_+)$, then $\widetilde F \in  \mathbb L^1(\mathbb R)$ and $\displaystyle \widehat{\widetilde F}(i.)$ is continuous.\\
Consequently if we establish~:
\begin{equation}
\label{majo-de-F-chap-tild}
 \left|\widehat{\widetilde F}(iq)\right|\leq\frac{C}{1+q^2}\eqp{,}
\end{equation}
then  \mbox{$ \widehat{\widetilde F}$}  will be  an element of $\mathbb L^1(\mathbb R)$.\\
It is proved in \cite{roynette-vallois-volpi-art1-2003} (Annex, Proposition A-2), that $\varphi(iq)=0$, $q\in \mathbb R$\; iff\;  $q=0$.\\
Therefore we are allowed to  replace $q$ by $iq$ in (\ref{trans-laplace-ruine}).\\
Using the identity~:
\begin{equation}
\int_{-1}^0(1+x)e^{-iqx}dx=- \frac{e^{iq}-1-iq}{q^2}\eqp{,}
\end{equation}
we then deduce~:
\begin{equation}
\label{F-tild-chapeau-bis}
\widehat {\widetilde F}(iq)
=\frac{1-e^{iq}}{q^2}+\frac{1}{\varphi(iq)}
\left[-c+\int_0^{+\infty}y\nu(dy)+\frac{i}{q}\int_{-k}^{0}(e^{-iqy}-1)\nu(dy)+RF(iq)\right]
\eqp{.}
\end{equation}
Since $\displaystyle \varphi(iq)\sim_{|q|\rightarrow+\infty}- \frac{q^2}{2}$, (\ref{majo-de-F-chap-tild}) is a consequence of the three inequalities~:
\begin{align}
&\left|\frac{1-e^{iq}}{q^2}\right|\leq\frac{2}{q^2}\eqp{,}\\
&\left|\frac{i}{q}\int_{-k}^0(e^{-iqy}-1)\nu(dy) \right|
\leq \int_{-k}^0|y|\nu(dy)\eqp{,}\\
& \left|RF(iq)\right|
\leq 2\int_{-k}^0 |y|\nu(dy)
\eqp{.}
\end{align}

\paragraph{{\bf Step 4}}The $n-2$ first derivatives of $\widehat{\widetilde F}$ belong to
$\mathbb L^1(\mathbb R)$.

\vskip6pt\noi
Obviously any $k$ derivative of $\displaystyle q\rightarrow \frac{1-e^{iq}}{q^2}$ is continuous,
 bounded by $\displaystyle \frac{C}{q^2}$, $|q| \geq1$, and consequently belongs
 to $\mathbb L^1(\mathbb R)$.\\
The second term in the right hand-side of (\ref{F-tild-chapeau-bis}) may be written
as  $\displaystyle \frac{N(q)}{\varphi(iq)}$.\\
Assumption (\ref{moment-p-fini}) implies that the $n$ first derivatives of $N$
are bounded, hence $\displaystyle \left|\frac{N(q)}{\varphi(iq)}\right|\leq\frac{C}{q^2}$, $|q|\geq 1$.\\
As for the asymptotic behaviour of $\displaystyle
\frac{N(q)}{\varphi(iq)}$ in a  neighborhood of $0$, it can be
proved by a similar reasoning  that this ratio is bounded, for any
$|q|\leq1$.

\paragraph{{\bf Step 5}} Proof of (\ref{majo-poly-de-F}).

\vskip6pt\noi Since the $n-2$ derivatives of   $q \rightarrow
\widehat {\widetilde F}(iq)$ belong to $\mathbb L^1(\mathbb R)$,
then

\begin{equation}
 x^{n-2} \widetilde F (x) =\frac{i ^{n-2}}{2\pi} \int_{\mathbb R}e^{iqx}\frac{d^{n-2}}{dq^{n-2}}\left(\widehat{\widetilde F}(iq)\right)dq\eqp{.}
\end{equation}
This identity directly implies (\ref{majo-poly-de-F}).
\finpreuve

%==========================================================================

\sectio{Normalized limit distribution of  $(T_x,K_x,L_x)$, as  \mbox{$x \rightarrow +\infty$}}
\label{sec:TCL}

%=========================================================================

In this section we investigate the limit behaviour of  $(T_x,K_x,L_x)$, as
 $x \rightarrow +\infty$. Recall that $K_x$ and $L_x$ are defined
 by (\ref{def-de-undershoot-overshoot}).  Since $x\rightarrow T_x$ is non-decreasing,
 $T_x$ needs to be normalized. We have to consider  three cases either  $\mathbb E(X_1)>0$,
 or $\mathbb E(X_1)<0$ or $\mathbb E(X_1)=0$.
Let us start with the case $\mathbb E(X_1)< 0$.

% ------------------THEOREME-------------------------
\begin{thm}

\label{thm:TCL-profit-net}
 \quad Under hypothesis of Theorem~\ref{thm:comp-asymp-de-F} and  $\mathbb E(X_1)< 0$, then,
conditionally on  \mbox{$\displaystyle \{T_x<+\infty\}$},
\mbox{$\displaystyle \left(\frac{1}{\sqrt{x}}
\left(T_x+\frac{x}{\varphi'(-\gamma_0(0))}\right),K_x,L_x\right)$}
converges in distribution to the 3-dimensional law

\noi $\displaystyle \mathcal
N\left(0;-\frac{\varphi''(-\gamma_0(0))}{\varphi'^3(-\gamma_0(0))}\right)
\otimes w$  \quad where $w$ is the probability measure on $\mathbb
R_+\times\mathbb R_+ $~:
\begin{align}
\label{loi-Kx-Lx}
 w(dk,dl)=&
\frac{-1}{\mathbb E(X_1)}\left[ \frac{\gamma_0(0)}{2} \; \delta_{0,0}(dk,dl)+ (e^{\gamma_0(0)l}-1)\un_{\{k\geq 0; l\geq 0\}}\;\nu_l(dk)\;dl\right.\nonumber\\[1ex]
& \left.
+\int_{\mathbb R_-}\left(\int_0^{-y}(1-e^{\gamma_0(0)(b+y)})F(0,0,0,b) n(b,dk,dl)db\;\un_{\{k\geq 0,l\geq 0\}}\right)\nu(dy)
\right]\eqp{,}
\end{align}

 \noi $\nu_l(dk)$ is the image of $\nu(dk)$ by the map $y\rightarrow
 y-l$,
 and $n(b,dk,dl)$ is the distribution of $(K_b,L_b)$
 conditionally on \mbox{$\displaystyle \{T_b<+\infty\}$}.\\
If moreover the support of $\nu$ is included in $[0;+\infty[$,
$w(dk,dl)$ is given explicitly :
\begin{equation}
\label{loi-K-L-vu+}
w(dk,dl)=\frac{-1}{\mathbb E(X_1)}\left[\frac{\gamma_0(0)}{2} \; \delta_{0,0}(dk,dl)+ (e^{\gamma_0(0)l}-1)\un_{\{k\geq0; l\geq 0\}}\;\nu_l(dk)\;dl
\right]
\end{equation}
\end{thm}

% --------------------------REMARQUE--------------------------------
\begin{rem}
\label{rem:TCL}
{\rm
\begin{enumerate}
%\item Result that $\gamma_0(0)=0$ and $T_x<+\infty$ a.s. if
%$\mathbb E(X_1)>0$.
\item  \mbox{$\displaystyle \mathcal N (0;\sigma^2)$} denotes the
Gaussian distribution with mean $0$ and variance $\sigma^2$. \item
A simular result has been proved by  Allan Gut (cf. \cite{Gut88},
page $102$, Theorem~$10.11$),
 where the L\'{e}vy process ($\X)$  is replaced by a random walk.
\item a) We observe that  time and positions become asymptotically independent. However
 the two components of the position are  not independent. We give a more
  complete description of the limit distribution of $(K_x,L_x)$
   in Proposition~\ref{prop:loi-L+K} below.\\
b) Obviously $K_x+L_x=X_{T_x}-X_{{T_x}^-}$ is the jump size of $(\X)$ at $T_x$. It is easy to deduce from Theorem~\ref{thm:TCL-profit-net} that
\mbox{$\displaystyle \left(\frac{1}{\sqrt{x}} \left(T_x+
\frac{x}{\varphi'(-\gamma_0(0))}\right), X_{T_x}-X_{{T_x}^-}\right)$} converges in distribution to \quad  \mbox{$\displaystyle \mathcal N \left(0;-\frac{\varphi''(-\gamma_0(0))}{\varphi'^3(-\gamma_0(0))}\right) \otimes w_0$}\quad where  $w_0$ is the probability measure on  $\mathbb R_+$~:
\begin{align}
 w_0(ds)=&
\frac{-1}{\mathbb E(X_1)}\left[ \frac{\gamma_0(0)}{2} \; \delta_{0}(ds)+
\frac{e^{\gamma_0(0)s}-1-\gamma_0(0)s}{\gamma_0(0)}\un_{\{s\geq 0\}}\;\nu(ds)\right.\nonumber\\[1ex]
& \left.
+\int_{\mathbb R_-}\left(\int_0^{-y}(1-e^{\gamma_0(0)(b+y)})F(0,0,0,b) n(b,ds)db\right)\nu(dy)\un_{\{s\geq 0\}}
\right]\eqp{,}
\end{align}
and $n(b,ds)$ is the distribution of $X_{T_b}-X_{{T_b}^-}$
conditionally on \mbox{$\displaystyle \{T_b<+\infty\}$}.\\
If the jumps of $(\X)$ are positive, $w_0(ds)$ can be simplified~:
\begin{equation}
w_0(ds)=\frac{-1}{\mathbb E(X_1)}\left[ \frac{\gamma_0(0)}{2} \;
\delta_{0}(ds)+
\frac{e^{\gamma_0(0)s}-1-\gamma_0(0)s}{\gamma_0(0)}\un_{\{s\geq
0\}}\;\nu(ds)\right].
\end{equation}
\item  a) If $\mathbb E(X_1)>0$, then $T_x<+\infty$ a.s. and the triplet
\mbox{$\displaystyle \left(\frac{1}{\sqrt{x}}
\left(T_x+\frac{x}{\varphi'(0)}\right),K_x,L_x\right)$} converges
 in distribution to  \mbox{$\displaystyle \mathcal N
 \left(0;-\frac{\varphi''(0)}{\varphi'^3(0)}\right) \otimes w$},
 where $w$ is defined by the relation obtained replacing $\gamma_0(0)$
   by $-\gamma_0^*(0)$ in
 (\ref{loi-Kx-Lx}).\\
 In particular if $(\X)$ has only positive jumps~:
\begin{equation}
w(dk,dl)=\frac{1}{\mathbb E(X_1)}\left[\frac{\gamma_0^*(0)}{2} \;
\delta_{0,0}(dk,dl)+ (1-e^{-\gamma_0^*(0)l})\un_{\{k\geq0; l\geq
0\}}\;\nu_l(dk)\;dl \right].
\end{equation}
b) In the third case~:  $\mathbb E(X_1)=0$, the normalization concerning $T_x$ has to be
  modified. In that case $\displaystyle \left(\frac{T_x}{x^2},K_x,L_x\right)$
   converges in distribution to $\varrho\otimes w$, where $\varrho$ denotes the law
   of the first hitting time of level $\displaystyle \sqrt{\frac{1}{\varphi''(0)}}$
   by a standard Brownian motion started at $0$.
   Recall that $\displaystyle \int_0^{+\infty} e^{-\theta x}\varrho(dx)=
   e^{-\sqrt{\frac{2\theta}{\varphi''(0)}}}$. Moreover
    $\displaystyle \varphi''(0)=1+\int_{\mathbb R}y^2\nu(dy)$).\\[1ex]
The probability measure $w$ on $\mathbb R_+\times\mathbb R_+$ is defined as follows~:
\begin{align}
\label{loi-Kx-Lx-egalite}
 w(dk,dl)=&
\frac{1}{\varphi''(0)}\left[ \vphantom{\int} \delta_{0,0}(dk,dl)+2l\un_{\{k\geq 0; l\geq 0\}}\;\nu_l(dk)\;dl\right.\nonumber\\[1ex]
& \left.
-2\int_{\mathbb R_-}\left(\int_0^{-y}(b+y) n(b,dk,dl)db\;\un_{\{k\geq 0,l\geq 0\}}\right)\nu(dy)
\right]\eqp{,}
\end{align}
where  $n(b,dk,dl)$ is the distribution of $(K_b,L_b)$.\\
This expression may be simplified if $(\X)$ has only positive jumps~:
\begin{equation}
w(dk,dl)=\frac{1}{\varphi''(0)}\left[ \delta_{0,0}(dk,dl)+2l\un_{\{k\geq 0; l\geq 0\}}\;\nu_k(dl)\;dk
\right]
\end{equation}
\end{enumerate}
}
\end{rem}

\noindent {\bf Proof of Theorem~\ref{thm:TCL-profit-net}}\mbox{}

\noindent Our approach is based on the following estimate :
\begin{equation}
\label{equival-de-F-bis}
F(\theta,\mu,\rho,x)\sim_{x\rightarrow+\infty} C_0(\theta,\mu,\rho)e^{-\gamma_0(\theta)x}\eqp{,}
\end{equation}
where $\displaystyle F(\theta,\mu,\rho,x)=\mathbb
E\left(e^{-\theta T_x-\mu K_x-\rho
L_x}\un_{\{T_x<+\infty\}}\right)$, and   $C_0(\theta,\mu,\rho)$
 (resp. $\gamma_0(\theta))$ is determined by
(\ref{residu}), (\ref{res-0-mu-cas-3}) (resp. (\ref{existence-de-kappa-bis})).\\
The three cases  $\mathbb E(X_1)>0$, $\mathbb E(X_1)<0$ and $\mathbb E(X_1)=0$ may be
treated with the same technic (cf. \cite{roynette-vallois-volpi-art1-2003}, Introduction, point 7.).
Therefore we only deal with  the third case~: $\mathbb E(X_1)=0$.\\
For simplicity we restrict ourselves to  the couple $(T_x, K_x)$,
 instead of the triplet $(T_x, K_x, L_x)$.\\
Since \mbox{$ \displaystyle \varphi(-\gamma_0(\theta))=\theta $}, $\varphi(0)=0$ and $\varphi'(0)=0$, taking the asymptotic expansion of $\varphi$ at $0$, of order $2$, we obtain~:
\begin{equation}
\theta =\varphi(-\gamma_0(\theta))=\varphi(-h)=
\frac{h^2}{2}\varphi''(0)+{\rm o}(h^2)
\eqp {.}
\end{equation}
Consequently~:
\begin{equation}
\gamma_0(\theta)=h= \sqrt{\frac{2\theta}{\varphi''(0)}}+{\rm o}(\sqrt{\theta})
\eqp{.}
\end{equation}
Recall that (\ref{equival-de-F-bis}) is uniform with respect
 to $\theta\in [0,\kappa]$. Then choosing $\rho=0$ and replacing $\theta$ by
   $\displaystyle \frac{\theta}{x^2}$ in (\ref{equival-de-F-bis}) bring  to~:
\begin{equation}
\mathbb E
\left(e^{-\theta \frac{T_x}{x^2}-\mu K_x}\right)
\sim_{x\rightarrow +\infty}
C_0(0,\mu)
e^{-\sqrt{\frac{2\theta}{\varphi''(0)}}}
\eqp{.}
\end{equation}
We have to check that $C_0(0,\mu)$ and $\displaystyle e^{-\sqrt{\frac{2\theta}{\varphi''(0)}}}$
are Laplace transforms of probability measures on $[0,+\infty[$. As for
 $\displaystyle e^{-\sqrt{\frac{1}{\varphi''(0)}}
\sqrt{2 \theta}}$, it is well known (cf. \cite{Karatszas91} page 96
formula (8.6)) that it is the Laplace transform of the first hitting time of level
$\displaystyle \sqrt{\frac{1}{\varphi''(0)}}$ by a standard Brownian motion started at $0$.\\
We modify the  identity (\ref{res-0-mu-cas-3}) (with $\rho=0$) via
\mbox{ $\displaystyle \frac{e^{-ay}-1+ay}{a^2}=-\int_0^y
(z-y)e^{-az}dz $}. Then    we obtain~:
\begin{align}
C_0(0,\mu)=&
\frac{1}{\varphi''(0)}
\left[1+ 2
\int_0^{+\infty}e^{-\mu z}\left(\int_{[z,+\infty[}(y-z)\nu(dy)\right)dz
\right.
\nonumber\\[1ex]
&\left.
-2\int_0^{+\infty}e^{-\mu z}\int_0^{+\infty} db \left(\int_{-\infty}^b (y+b)e^{-(y+b)x}\nu(dy)\right) n(b,dz)  \right]\nonumber\\
&=\int_0^{+\infty} e^{-\mu z} w(dz)
\eqp{.}
\end{align}
\finpreuve
%

% -----------------REMARQUE--------------
\begin{rem} \label{rasymp3} We would like to point out that
the asymptotic development of $F(\theta , \mu, \rho,x)$ given by
(\ref{asymp4})  gives the rate of convergence of $\displaystyle
\left(\frac{1}{\sqrt{x}}
\left(T_x+\frac{x}{\varphi'(-\gamma_0(0))}\right)\right)$ to the
Gaussian distribution. Suppose that $\mathbb{E}[X_1]>0$. Let $A_x$
be the distribution function of $\displaystyle
\left(\frac{1}{\sqrt{x}}
\left(T_x+\frac{x}{\varphi'(-\gamma_0(0))}\right)\right)$ and
$\widehat{A}_x$ be its characteristic function :
$$
\begin{array}{ccll}
  A_x(t) & = &\displaystyle  \mathbb{P} \left(\frac{1}{\sqrt{x}}
\big(T_x+\frac{x}{\varphi'(-\gamma_0(0))}\big) \leq t\right),& t \in \mathbb{R} \\
&&&\\
  \widehat{A}_x(\theta) & =\displaystyle
  & \mathbb{E}\Big[e^{ \frac{i\theta}{\sqrt{x}}(T_x+\frac{x}{\varphi'(-\gamma_0(0))})}
  \Big],
  & \theta \in \mathbb{R}. \\
\end{array}
$$
\noi Thanks to (\ref{asymp4}), it is no difficult to check that,
if $x$ is large enough :
\begin{equation}\label{conv1}
    \sup_{\theta >0}|\widehat{A}_x(\theta)-\widehat{A}(\theta)|
    \leq \frac{k\theta}{\sqrt{x}},
\end{equation}
\noi where  $k$ is a positive constant and $\widehat{A}$ is the
characteristic function of $\displaystyle \mathcal
N\left(0;-\frac{\varphi''(-\gamma_0(0))}{\varphi'^3(-\gamma_0(0))}\right)$.

\noi Berry-Essen's inequalities (see for instance \cite{Loeve}, p
285) implies the existence of two positives constants $k_1$ and
$k_2$ such that :
\begin{equation}\label{conv2}
\sup_{t \in \mathbb{R}}|A_x(t)-A(t)|\leq k_1\int_0^a
\frac{k\theta}{\sqrt{x}}\frac{d\theta}{\theta} +\frac{k_2}{a},
\quad \mbox{for any } \ a>0,
\end{equation}
\noi where $A$ is the distribution function of the  previous
Gaussian distribution.

\noi Choosing $a=x^{1/4}$, we obtain :
\begin{equation}\label{conv3}
    \sup_{t \in \mathbb{R}}|A_x(t)-A(t)|\leq \frac{k_3}{x^{1/4}},
    \ x\geq 1,
\end{equation}
for some $k_3>0$.
\end{rem}
% ---------------------------------

\noi We would like to provide a stochastic interpretation of the probability
measure $w$ defined by (\ref{loi-K-L-vu+}). Let $(K,L)$ be a two-dimensional r.v.\
 with probability distribution $w$. Obviously $\{K=0\}=\{L=0\}$ and this event
 occurs with probability $\displaystyle \frac{-\gamma_0(0)}{2\mathbb E(X_1)}$. Conditionally on
  $\{L>0\}$, the distribution of $(K,L)$ is of type
  $\displaystyle \alpha\left(e^{\gamma_0(0)l}-1\right)\un_{\{k>0;l>0\}}\nu_l(dk) dl$
   where $\nu_l$ is the positive measure defined in Theorem~\ref{thm:TCL-profit-net}.\\
This leads us to consider the positive measure~:
\begin{equation}
\label{def-w*}
w^*(dk,dl)=\alpha(e^{\gamma l}-1)\un_{\{k>0;l>0\}}\nu_l(dk)dl\eqp{,}
\end{equation}
where $\gamma >0$, ${\nu}_l$ is the image of $\nu$ by
$y\rightarrow y-l$ and $\nu$ is a positive measure  on
$]0;+\infty[$ satisfying ~:
\begin{equation}
\int_0^{+\infty}(e^{\gamma k}-1-\gamma k)\nu(dk) < +\infty \eqp{,}
\end{equation}
 $\alpha$ being the normalization factor :
 $ \displaystyle \alpha=\frac{\gamma}{\displaystyle \int_0^{+\infty}(e^{\gamma k}-1-\gamma k)\nu(dk)}$.

\begin{prop}
\label{prop:loi-L+K}
Let $(K^*,L^*)$ be a two dimensional r.v.\  with distribution $w^*$ defined by (\ref{def-w*}) . Then $L^*$ has a density function given by $\displaystyle \alpha (e^{\gamma l}-1)\nu([l,+\infty[)\un_{\{l>0\}}$. Conditionally to $L^*=l$, the distribution of $S^*=L^*+K^*$ is $\displaystyle \frac{1}{\nu([l,+\infty[)}\un_{\{s>l\}}\nu(ds)$.
\end{prop}
We have considered  $w$ defined by (\ref{loi-Kx-Lx}), but  a
similar analysis can be developed with $w$ satisfying
(\ref{loi-Kx-Lx-egalite}).

\clearpage

\def\refname{References}
\bibliographystyle{plain}
\bibliography{BiblioLevy}

\begin{thebibliography}{10}

\bibitem{BertoinDoney94}
J.~Bertoin and R.A. Doney.
\newblock Cram\'er's estimate for {L}\'{e}vy processes.
\newblock {\em Statistics and Probability Letters}, 21:363--365, 1994.

\bibitem{Doney91}
R.A. Doney.
\newblock Hitting probabilities for spectrally positive {L}\'{e}vy processes.
\newblock {\em J.London Math. Soc.}, 44:566--576, 1991.

\bibitem{DufresneGerber90}
F.~Dufresne and H.U. Gerber.
\newblock Risk theory for the compound {P}oisson process that is perturbed by
  diffusion.
\newblock {\em Insurance: Mathematics and Economics}, 10(1991):51--59, 1990.

\bibitem{FrolovaKabanov01}
A.G. Frolova, Y.M. Kabanov, and S.M. Pergamenshchikov.
\newblock In the insurance business risky investments are dangerous.
\newblock {\em Preprint, Univ. de Franche Comt\'{e}, Besan\c{c}on}, 2001.

\bibitem{FurrerSchmidli94}
H.J. Furrer and H.~Schmidli.
\newblock Exponential inequalities for ruin probabilities of risk processes
  perturbed by diffusion.
\newblock {\em Insurance: Mathematics and Economics}, 15:23--36, 1994.

\bibitem{GerberLandry98}
H.U. Gerber and B.~Landry.
\newblock On the discounted penalty at ruin in a jump-diffusion and the
  perpetual put option.
\newblock {\em Insurance: Mathematics and Economics}, 22:263--276, 1998.

\bibitem{Gut88}
A.~Gut.
\newblock {\em Stopped Random Walks. Limit Theorems and Applications}, volume~5
  of {\em Applied Probability}.
\newblock Springer-Verlag, New York, 1988.

\bibitem{Karatszas91}
I.~Karatszas and S.E. Shreve.
\newblock {\em Brownian motion and Stochastic Calculus}.
\newblock Springer-Verlag, New York, second edition, 1991.
\newblock Graduate Texts in Mathematics, 113.

\bibitem{KluppelbergKyprianou03}
C.~Kluppelberg, A.E. Kyprianou, and R.A. Maller.
\newblock {R}uin {P}robabilities and {O}vershoot for {G}eneral {L}{\'{e}}vy
  {I}nsurance {R}isk {P}rocesses.
\newblock preprint, 2003.

\bibitem{Loeve}
M.~Lo{\`e}ve.
\newblock {\em Probability theory}.
\newblock Third edition. D. Van Nostrand Co., Inc., Princeton, N.J.-Toronto,
  Ont.-London, 1963.

\bibitem{Paulsen01}
J.~Paulsen.
\newblock On {C}ram\`{e}r-like asymptotics for risk processes with stochastic
  return on investments.
\newblock {\em Preprint}, 2001.

\bibitem{roynette-vallois-volpi-art1-2003}
B.~Roynette, P.~Vallois, and A.~Volpi.
\newblock {L}\'{e}vy processes~: Hitting time, overshoot and undershoot; {I}-
  functional equations.
\newblock {\em Preprint}, 2003.

\bibitem{Schmidli95}
H.~Schmidli.
\newblock Cram\`{e}r-{L}undberg approximations for ruin probabilities of risk
  processes perturbed by diffusion.
\newblock {\em Insurance: Mathematics and Economics}, 16:135--149, 1995.

\bibitem{Schmidli01}
H.~Schmidli.
\newblock Distribution of the first ladder height of a stationary risk process
  perturbed by $\alpha$-stable {L}\'{e}vy motion.
\newblock {\em Insurance: Mathematics and Economics}, 28:13--20, 2001.

\end{thebibliography}

\end{document}